\documentclass[a4paper,10pt]{article}
\pdfoutput=1
\usepackage[french, english]{babel}
\usepackage[utf8]{inputenc}
\usepackage[T1]{fontenc}
\usepackage{amsmath,amssymb,amsthm,mathtools}
\usepackage{xpatch}
\usepackage[backend=biber,doi=false,url=false,isbn=false,maxbibnames=99]{biblatex}
\usepackage{geometry}
\usepackage{IEEEtrantools}
\usepackage{tikz}
\usetikzlibrary{math,calc,patterns,decorations.markings}
\usepackage{pgfplots}
\pgfplotsset{compat=1.18}
\usepackage[titletoc,title]{appendix}
\usepackage{stix2}
\usepackage{dsfont}
\usepackage{microtype}
\usepackage[font=footnotesize, margin=1em]{caption}
\usepackage{csquotes}
\usepackage{hyperref}
\usepackage{cleveref}

\DeclareFieldFormat{eprint:hal}{%
  \href{https://hal.archives-ouvertes.fr/#1}{%
  \nolinkurl{#1}}}
\AtEveryBibitem{\clearfield{day}}
\AtEveryBibitem{%
  \ifentrytype{book}
    {\clearfield{pagetotal}}
    {}%
}
\definecolor{sienne}{RGB}{136, 45, 23}
\hypersetup{
colorlinks = true,
linkcolor = sienne
}
\allowdisplaybreaks[1]

\newtheorem{theorem}{Theorem}[section]
\newtheorem{prop}[theorem]{Proposition}
\newtheorem{corollary}[theorem]{Corollary}
\newtheorem{lemma}[theorem]{Lemma}
\theoremstyle{definition}
\newtheorem{definition}[theorem]{Definition}
\theoremstyle{remark}

\makeatletter
\patchcmd{\@IEEEeqnarray}{\relax}{\relax\intertext@}{}{}
\makeatother

\newcounter{proofstep}[theorem]
\newcommand{\step}[1]{%
\refstepcounter{proofstep}%
\vskip-\lastskip\medskip\noindent\textit{Step \arabic{proofstep}: #1. ---}}
\xapptocmd\proof{\setcounter{proofstep}{0}}{}{}

\newenvironment{delayedproof}[1]
 {\begin{proof}[\raisedtarget{#1}Proof of \cref{#1}]}
 {\end{proof}}
\newcommand{\raisedtarget}[1]{%
  \raisebox{\fontcharht\font`P}[0pt][0pt]{\hypertarget{#1}{}}%
}
\newcommand{\proofref}[1]{\hyperlink{#1}{proof}}

\DeclareMathOperator{\distance}{d}

\DeclareMathOperator{\range}{Range}
\DeclareMathOperator{\supp}{supp}

\DeclareMathOperator{\Rg}{Rg}

\DeclareMathOperator{\Li}{Li_2}
\newcommand{\dz}[1][z]{\partial_{#1}}
\newcommand{\dzb}[1][z]{\partial_{\overline{#1}}}
\DeclareMathOperator{\Real}{Re}

\newcommand{\bigO}{O}

\newcommand{\diff}[1][-3]{\mathop{}\mkern#1mu\mathrm{d}}
\newcommand{\set}{\mathbb}
\newcommand{\N}{{\set{N}}}
\newcommand{\Z}{{\set{Z}}}
\newcommand{\R}{{\set{R}}}
\newcommand{\C}{{\set{C}}}
\newcommand{\T}{{\set{T}}}
\newcommand{\D}{{\set{D}}}
\newcommand{\Dext}{{\set{D}_{\mathrm{e}}}}

\newcommand{\Cinf}{{\widehat{\C}}}
\newcommand{\Deinf}{{\widehat{\D}_{\mathrm{e}}}}
\newcommand{\Hol}{H}

\newcommand{\cl}[1]{{#1^{\mathsf{cl}}}}

\newcommand{\an}{{\mathrm a}}

\newcommand{\FT}[1][T]{{\mathcal{F}_{#1}}}
\newcommand{\FTp}[1][T]{{\mathcal{F}_{#1}^+}}

\newcommand{\NC}{{\mathcal{NC}}}
\newcommand{\NCL}{{\NC_{L^2}}}
\newcommand{\NCH}{{\NC_{H^2}}}
\newcommand{\Reach}{{\mathcal{R}}}
\newcommand{\ReachH}{{\Reach_{H^2}}}
\newcommand{\ReachL}{{\Reach_{L^2}}}

\newcommand{\Dirichlet}{\mathcal{D}_2}

\newcommand{\coloneqq}{\mathrel{\mathord{:}\mathord=}}
\newcommand{\eq}{\Leftrightarrow}

\newcommand{\eu}{\mathrm e}

\newcommand{\iu}{\mathrm i}

\definecolor{omegac}{rgb}{0.5,0,0}
\definecolor{Omegac}{rgb}{1,0.2,0.2}
\definecolor{Omegac2}{rgb}{1,0.5,0.5}
\definecolor{diskc}{rgb}{1,1,0}

\title{Control of the half-heat equation}
\author{Andreas Hartmann\thanks{Institut de Mathématiques Bordeaux, Université de Bordeaux, Bordeaux, France. \texttt{andreas.hartmann@math.u-bordeaux.fr}.}, Armand Koenig\thanks{Institut de Mathématiques Bordeaux, Université de Bordeaux, Bordeaux, France. \texttt{armand.koenig@math.u-bordeaux.fr}.}}

\begin{document}
\maketitle

\begin{abstract}
    In this paper we investigate null-controllable initial states of the half heat equation controlled from a sub-arc $\omega$ of the unit circle. We also study the projection on positive frequencies of the half-heat equation. For this projected half-heat equation, we obtain necessary as well as sufficient conditions for an initial condition to be null-controllable. These conditions, which are almost sharp, are expressed in term of projections on positive frequencies of functions supported on $\omega$. From these results, and with the help of classical results on sum of holomorphic and anti-holomorphic functions, we also treat the (unprojected) half-heat equation. Surprisingly, without using our conditions on null-controllable states, we are able to show that the space of null-controllable functions does not depend on time by using a result of separation of singularities for holomorphic functions.
\end{abstract}

\section{Introduction}
\subsection{Setting}
Let $\T$ be the unit complex circle. Functions in $L^2(\T)$ can be written as $f(\eu^{\iu x}) = \sum_{n\in\Z} \widehat f(n) \eu^{\iu nx}$, where $(\widehat f(n))_n \in \ell^2(\Z)$. With this notation, we will frequently identify $\T$ and $\R/2\pi\Z\sim [0,2\pi)$, and use in a similar fashion the words arc and interval for a connected component of $\T$. Let $|D|$ be the unbounded operator on $L^2(\T)$ defined by
\begin{equation}\label{eq-def-D}
    |D|f(\eu^{\iu x}) = \sum_{n\in \Z}|n| \widehat f(n) \eu^{\iu n x}.
\end{equation}
Given an open strict interval $\omega\subset \T$, we consider the control system (that we will call ``$L^2$ control system''):
\begin{equation}\label{eq-half-L2}
    \left\{\begin{IEEEeqnarraybox}[][c]{l"l}
        (\partial_t + |D|)f(t,\eu^{\iu x}) = \mathds 1_\omega u(t,\eu^{\iu x}),& t\geq 0,\ \eu^{\iu x}\in \T.\\
        f(0,\eu^{\iu x}) = f_0(\eu^{\iu x}),& \eu^{\iu x}\in\T,\ f_0  \in L^2(\T).
    \end{IEEEeqnarraybox}\right.
\end{equation}
As an intermediary step, we also consider the following ``$H^2$ control system'', where $P_+ f(\eu^{\iu x}) = \sum_{n\geq0} \widehat f(n)\eu^{\iu n x}$ is the Riesz projector and $H^2 = P_+ L^2(\T)$:
\begin{equation}\label{eq-half-H2}
    \left\{\begin{IEEEeqnarraybox}[][c]{l"l}
        (\partial_t + |D|)f(t,\eu^{\iu x}) = P_+\mathds 1_\omega u(t,\eu^{\iu x}),& t\geq 0,\ \eu^{\iu x}\in \T.\\
        f(0,\eu^{\iu x}) = f_0(\eu^{\iu x}),& \eu^{\iu x}\in\T,\ f_0 \in H^2.
    \end{IEEEeqnarraybox}\right.
\end{equation}
These are control systems with state $f(t,\cdot)\in L^2(\T)$ (respectively $H^2$) and control $u\in L^2([0,T]\times \omega)$. We are interested in which initial states can be steered to $0$, namely:
\begin{definition}\label{def-NC}
    We say that $f_0\in L^2(\T)$ (resp.\ in $H^2(\T)$) is null-controllable (or steerable to $0$) from $\omega$ in time $T$ in $L^2$ (resp.\ in $H^2$) if there exists $u\in L^2([0,T]\times \omega)$ such that the solution $f$ of the $L^2$ control system~\eqref{eq-half-L2} (resp.\ the $H^2$ control system~\eqref{eq-half-H2}) with initial condition $f_0$ is such that $f(T,\cdot) = 0$.

    We denote by $\NCL(\omega,T)$ (resp.\ $\NCH(\omega,T)$) the set of initial conditions that are steerable to $0$ in $L^2$ (resp.\ $H^2$) from $\omega$ in time $T$.
\end{definition}

Using Fourier series, we can check that the operator $|D|$ with domain $W^{1,2}(\T)$ is self-adjoint maximal monotone on $L^2(\T)$, hence it generates an analytic semigroup $S(t)$ on $L^2(\T)$. If $f_0(x)=\sum_{n\in\Z} \widehat{f_0}(n) \eu^{\iu nx}$, this semi-group is given by
\begin{equation}\label{eq-semigroup}
    S(t)f_0(\eu^{\iu x})=\sum_{n\in\Z} \widehat{f_0}(n) \eu^{-|n|t}\eu^{\iu nx}.
\end{equation}
In other words, the Cauchy problems \eqref{eq-half-L2} and \eqref{eq-half-H2} are well-posed.

\subsection{Overview of the results}

% \subsubsection{The \texorpdfstring{$H^2$}{}-control system}

The main idea of this article is that free solutions of the half-heat equation are sums of a holomorphic function and an anti-holomorphic function in $\eu^{-t+\iu x}$, as can be seen in the formula for the semigroup~\eqref{eq-semigroup}.

We start by discussing the null-controllable states for the $H^2$-control system~\eqref{eq-half-H2}, where free solutions are now holomorphic in $\eu^{-t+\iu x}$.

\begin{figure}
    \begin{minipage}[c]{0.4\textwidth}
        \begin{center}
            \tikzmath{\Rad=1.7; \rad=1/\Rad; \a=10; \b=60; \abmid = (\a+\b)/2; \Rmid = (2+3*\Rad)/5; \rmid = (1+\rad)/2; }
\begin{tikzpicture}[scale=2, every node/.style={fill=white, fill opacity = 0.3, inner sep=1pt, text=black, text opacity=1}]
    \draw (0,0) circle[radius=1];
    
    \draw[fill=Omegac, fill opacity = 0.5] (\a:1) -- (\a:\Rad)
        arc[start angle=\a, end angle=\b, radius=\Rad]
        -- (\b:1)
        arc[start angle=\b, end angle=\a, radius=1];
        
    \draw[fill=diskc, fill opacity=0.5] (\a:1) -- (\a:\rad)
        arc[start angle=\a, end angle=\b, radius=\rad]
        -- (\b:1)
        arc[start angle=\b, end angle=\a, radius=1];
    
    \draw[|-|, ultra thick, omegac] (\a:1) 
    arc[start angle=\a, end angle = \b, radius = 1] 
    node[pos = 0.5, above right]{$\omega$};
    
    \node at (\abmid:\Rmid) {$\Omega_T$};
    \node at (\abmid:\rmid) {$\Omega_T^*$};
    \node at (135:0.5){$\D$};
    \node[below right] at (-45:1) {$\T$};
    
    \draw[->] (-1.2,0) -- (1.2,0);
    \draw[->] (0,-1.2) -- (0,1.2);
\end{tikzpicture}
        \end{center}
    \end{minipage}\hfill%
    \begin{minipage}[c]{0.6\textwidth}
        \caption{Illustration of $\Omega_T^*$ (yellow partial ring inside of $\D$) and $\Omega_T$ (red partial ring outside $\D$).}
        \label{fig-Omega}
    \end{minipage}
\end{figure}
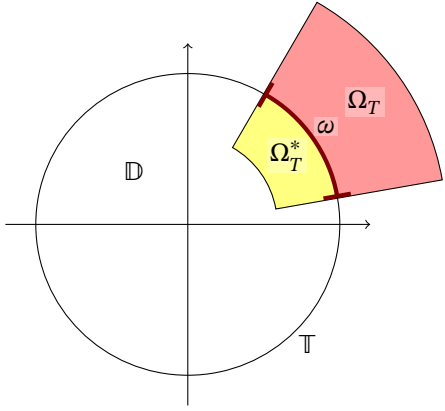

The starting point of our analysis is the usual equivalence between null-controllability and the so-called observability inequality~\cites[Theorem 2.44]{Coron07}[Theorem 11.2.1]{TW09}. A central role in our analysis will be played by Bergman spaces. We recall that for a domain $\Omega\subset \C$, the Bergman space $A^2(\Omega)$ is the set of functions $f$ holomorphic on $\Omega$ such that $\|f\|_{A^2(\Omega)}^2=\int_{\Omega}|f|^2\diff A<+\infty$, where $A$ is planar Lebesgue measure (weighted versions are defined integrating against the corresponding weight).  We refer to, e.g., \cite[\S18.1]{Conway95} for an elementary introduction to Bergman spaces. Combined with the change of variables $z = \eu^{-t+\iu x}$, we are able prove in \cref{sec-obs} the following characterization of null-controllable states:
\begin{prop}\label{th-MainEstim}
    Let $\omega$ be a strict open subset of $\set T$ and $T>0$. Let $\Omega_T \coloneqq \{ z\in \set C,\ 1<|z|<\eu^T,\ z/|z| \in \omega\}$ (see \Cref{fig-Omega}). Let $f_0 \in H^2$. The following assertions are equivalent:
    \begin{itemize}
        \item $f_0\in \NCH(\omega,T)$
        \item there exists $C_{f_0,T}>0$ such that for every polynomial $p\in \C[X]$,
              \begin{equation}\label{MainEstim}
                  |\langle f_0,p\rangle_{H^2}| \leq C_{f_0,T}\|p\|_{A^2(\Omega_T)}.
              \end{equation}
    \end{itemize}
    If these assertions are satisfied, we can choose the control $u$ that steers $f_0$ to $0$ such that $\|u\|_{L^2([0,T]\times\omega)} \leq  \eu^T C_{f_0,T}$, where $C_{f_0,T}$ is the constant that appears in the second assertion.
\end{prop}

Our first main result states that $\NCH(\omega,T)$ does not depend on $T$:
\begin{theorem}\label{th-NCH-T}
    Let $\omega$ be an interval of $\T$. For every $T,T'>0$, $\NCH(\omega,T) = \NCH(\omega,T')$.
\end{theorem}
The \proofref{th-NCH-T} of \cref{th-NCH-T} relies on a theory called \emph{separation of singularities}, that we recall in \cref{th-SepSing-est}. In view of \cref{th-NCH-T}, from now on, we will omit the $T$ in $\NCH(\omega,T)$ and write $\NCH(\omega)$.

We will actually prove a more precise result in \cref{th-NCH-T-2} which gives an estimate of the cost of the control: given $f\in \NCH(\omega)$, there is $C$ such that the function $u$ steering the system to zero in time $T$ satisfies
\[
 \|u\|_{L^2([0,T]\times\omega)} \leq C\max(1,T^{-3/2}).
\]

According to our definition of $H^2(\T) = P_+L^2(\T)$, functions in $H^2(\T)$ have the form $f\colon\eu^{\iu x}\in \T \mapsto \sum_{n\geq 0} \widehat f(n) \eu^{\iu n x}$. Since only non negative $n$ are present here, we can actually define $f(z) = \sum_{n\geq 0} \widehat f(n) z^n$ for $z\in\D \coloneqq \{z\in \C,\ |z|<1\}$. This way, we recover the usual definition of the Hardy space $H^2$ (see, e.g., \cite[Definition~17.7 \&{} Theorem~ 17.12]{Rudin86}), with the usual identification of $f\in H^2$ (seen as holomorphic in $\D$) and of its radial limits $f(\eu^{\iu \theta}) \coloneqq \lim_{r\to1^-} f(r\eu^{\iu \theta})$~\cite[Theorem~17.11]{Rudin86}. We will make this identification throughout the article.

If $E\subset \C$, we denote the closure of $E$ by $\cl E$. Our main necessary condition for a function $f_0\in H^2$ to be null-controllable is the following:
\begin{theorem}\label{th-NCH-holom-necessary}
    We denote the Riemann sphere by $\Cinf$. If $f_0 \in \NCH(\omega)$, then $f_0$ can be extended holomorphically to $\Cinf\setminus \cl\omega$. Denoting this extension also by $f_0$, we have
    \begin{enumerate}
        \item $f_0(\infty) =   0$, %f_0'(\infty) =
        \item $\int_{|z|>1} |f_0'(z)|^2 \diff A(z) <+\infty$.
    \end{enumerate}
\end{theorem}
The \proofref{th-NCH-holom-necessary} of \cref{th-NCH-holom-necessary} mostly consists in testing the observability inequality of \cref{th-MainEstim} on the \emph{Hardy reproducing kernel} (see definition below).
We will call the space of holomorphic functions on $\Dext$ satisfying (i) and (ii) the Dirichlet space of the extended outer disk and denote it by $\Dirichlet(\Dext)$. 
%In other words, a null-controllable initial state is in $H^2$ in the unit disk and in $\Dirichlet(\Dext)$ ($\Dext$ is the ``exterior of the unit disk'', i.e., $\{z\in\C,\ |z|>1\}$). 
Using a natural connection with holomorphic functions on $\Cinf\setminus \T$~\cite[Theorem~5.3.1]{CMR06}, we will see that every $f_0\in\NCH(\omega)$ can be written as the so-called \emph{Cauchy-transform} of a function supported on $\omega$. More precisely given $g\in L^2(\omega)$, we can write:
\begin{equation}\label{CauchyTrans}
 \int_{\omega}\frac{g(w)}{1-\overline{w}z}\lvert\diff w|=
 \begin{cases}
    P_+g(z)=\sum_{n\ge 0}\widehat{g}(n)z^n, & |z|<1,\\
    -P_-g(z)=-\sum_{n>0}\dfrac{\widehat{g}(-n)}{z^n}, & |z|>1,
 \end{cases}
\end{equation}
and the left-hand side is called the \emph{Cauchy transform of $g$} (here $P_-=I-P_+$ with its natural extension to $\D_e$). With this, we are able to rephrase \cref{th-NCH-holom-necessary} as:
\begin{theorem}\label{th-CondNec}
    Let $f_0 \in \NCH(\omega)$. There exists $g\in L^2(\T)$, such that
    \begin{itemize}
        \item $\supp g\subset \cl\omega$
        \item $\sum\limits_{n<0} |n| |\widehat g(n)|^2 < +\infty$
        \item $f_0 = P_+ g$.
    \end{itemize}
\end{theorem}

Another consequence of \cref{th-NCH-holom-necessary} is that we cannot steer non zero analytical initial states to zero.
\begin{corollary}\label{th-analy-non-nc}
    Let $\omega$ be a strict open subset of $\T$. If $f_0\in H^2(\T)$ extends holomorphically to a disk larger than $\D$ and $f_0\in \NCH(\omega)$, then $f_0 = 0$.
\end{corollary}

We also get a sufficient condition for an initial condition to be null-controllable:
\begin{theorem}\label{th-CondSuff}
    Let $\omega$ be an interval of $\T$ with $\partial\omega=\{\zeta_1,\zeta_2\}$. Suppose $g\in L^2(\T)$ with
    \begin{itemize}
        \item $\supp g \subset \cl\omega$
        \item $\sum\limits_{n<0} |n| |\widehat g(n)|^2 <+\infty$
        \item $\int_{\T}\frac{|g(\eu^{\iu t})|^2}{|(\eu^{\iu t}-\zeta_1)(\eu^{\iu t}-\zeta_2)|}\diff t<\infty.$
    \end{itemize}
    Then $P_+g\in \NCH(\omega)$.
\end{theorem}

As a matter of fact, the \proofref{th-CondSuff} will rely on a necessary and sufficient condition expressed in terms of pseudo-Carleson measures (\cref{CNS}). We refer to Xiao's article~\cite{Xiao00} for a discussion on pseudo-Carleson measure.

In the following, $W^{k,p}$ are the usual Sobolev spaces, see, e.g.,~\cite[Definition~7.14]{Leoni17} for $k\in \N$ and \cite[Definition~1.12]{Leoni23} for $0<k<1$. Note that according to \cref{th-analy-non-nc}, the complementary set of $\NCH(\omega)$ is dense in every $H^2(\T)\cap W^{k,2}(\T)$. We can also deduce from \cref{th-CondSuff} that $\NCH(\omega)$ is dense in every $H^2\cap W^{k,2}(\T)$.
\begin{corollary}\label{th-dense}
    Let $\omega$ be a non-empty open interval of $\T$.
    For every $k\in \N$, the set $\NCH(\omega)\cap W^{k,2}(\T)$ is dense in $H^2\cap W^{k,2}(\T)$.
\end{corollary}
The \proofref{th-dense} of \cref{th-dense} is found in \cref{sec-h2}. Note that $W^{k,2}(\T)\cap H^2(\T)$ is a closed subspace of $W^{k,2}(\T)$, and we endow $W^{k,2}(\T)\cap H^2(\T)$ with the $W^{k,2}(\T)$ scalar product.

There is still a gap between the necessary condition of \cref{th-CondNec} and the sufficient condition of \cref{th-CondSuff}. However, if we assume a priori regularity on $f_0$, we get a necessary and sufficient condition:
\begin{corollary}\label{th-CNS-W12}
    Let $f_0\in H^2\cap W^{1/2,2}(\T)$. Then $f_0\in \NCH(\omega)$ if and only if there exists $g\in L^2(\T)$ such that
    \begin{itemize}
        \item $\supp g \subset \cl\omega$
        \item $g\in W^{1/2,2}(\T)$
        \item $f_0 = P_+g$.
    \end{itemize}
\end{corollary}
We should emphasize that in \cref{th-CondSuff} the Sobolev type condition involves only the negative Fourier coefficients (corresponding to $\Dirichlet(\Dext)$) while in this corollary it concerns both positive and negative coefficients. Observe also that the first two conditions on $g$ in this corollary can be synthesized as $g\in W^{1/2,2}_{00}(\omega)$. Indeed, $W^{1/2,2}_{00}(\omega)$ is by definition the set of functions in $L^2(\omega)$ such that their extension by $0$ on $\T$ is in $W^{1/2,2}(\T)$~\cite[Definition 1.85]{Leoni23}. 
The \proofref{th-CNS-W12} of \cref{th-CNS-W12} mostly consists in plugging \cref{th-CondNec,th-CondSuff} with some Hardy-type inequality satisfied by functions in $W^{1/2,2}_{00}$~\cite[Theorem~1.87]{Leoni23}.

We can also get another equivalent condition to $f_0\in\NCH(\omega)$ although less tractable than the corollary above (it involves harmonic extension to $\Omega_T$). Using Stokes formula, we will see in \cref{sec-h2-poisson} that the following condition seems straightforward (but the rigorous proof takes quite the roundabout way):
\begin{theorem}\label{th-CNS}
    Let $g \in L^2(\omega)$, and let $u$ be the solution of $\Delta u = 0$ on $\Omega_T$, $u(z)=zg(z)$ on $\omega$, $u=0$ on $\partial \Omega_T\setminus\omega$ (see \cref{sec-conformal-and-poisson,th-poisson-lipschitz} for the definition). The following assertions are equivalent:
    \begin{itemize}
        \item $P_+g\in\NCH(\omega)$,
        \item $\dz u \in L^2(\Omega_T)$.
    \end{itemize}
\end{theorem}

The study of the $H^2$ control system~\eqref{eq-half-H2} is the stepping stone to the study of the full $L^2$ control system~\eqref{eq-half-L2}. We are able to prove in \cref{sec-l2} the following link between $\NCH(\omega)$ and $\NCL(\omega,T)$:
\begin{theorem}\label{th-NCH-sub-NCL}
    Let $\omega$ be a strict interval of $\T$ and let $T>0$. Then, for every $f_0\in L^2(\T)$.
    \[f_0 \in\NCL(\omega,T) \eq \left\{\begin{array}{l}
            P_+ f_0 \in \NCH(\omega) \\
            P_+ \overline{f_0} \in \NCH(\omega).
        \end{array}\right.
    \]
\end{theorem}
Recall that free solutions of the half-heat equation are sums of holomorphic and anti-holomorphic functions in $\eu^{-t+\iu x}$. If we are to translate the results on $\NCH(\omega)$ into results on $\NCL(\omega,T)$, we need a way to properly ``separate'' the dynamics of the positive and negative frequencies just by looking at the control region $\omega$. Fortunately, an inequality by Friedrichs~\cite{Friedrichs37} turns out to be just the right tool, and is the main technical ingredient of the \proofref{th-NCH-sub-NCL} of \cref{th-NCH-sub-NCL}.

Theorem \ref{th-NCH-sub-NCL} allows us to study the set $\NCL(\omega,T)$ by using the results on $\NCH(\omega,T)$. First, since $\NCH(\omega,T)$ does not depend on $T$ (\cref{th-NCH-T}), the right-hand side in the previous equivalence does not depend on $T$. Hence:
\begin{corollary}\label{th-NCL-T}
    Let $\omega$ be an interval of $\T$. The set $\NCL(\omega,T)$ does not depend on $T$: for all $T',T > 0$, $\NCL(\omega,T) = \NCL(\omega,T')$.
\end{corollary}

We will also prove the following analog of \cref{th-analy-non-nc}.
\begin{corollary}\label{th-NCL-non-NC}
    Let $\omega$ be a strict interval of $\T$ and let $T>0$. If $f_0\in L^2(\T)$ is analytic on $\T$ and $f_0\in \NCL(\omega,T)$, then $f_0 = 0$.
\end{corollary}

\Cref{th-CNS-W12} does not transfer so nicely to the $L^2$ control system (we would find that $\NCL(\omega)\cap W^{1/2,2}(\T) = \{P_+ g_+ + (1-P_+)g_-,\ g_\pm \in W^{1/2,2}_{00}(\T),\ \int_\T g_+ = \int_\T g_- \}$; we leave the proof to the interested reader), but we can still prove that the space of null-controllable initial data is dense:
\begin{corollary}\label{th-NCL-dense}
    For every $k\in \N$, the set $\NCL(\omega)\cap W^{k,2}(\T)$ is dense in $W^{k,2}(\T)$.
\end{corollary}

\subsection{Discussion and examples} 
Let us give some examples of functions that are or are not null-controllable.

Our necessary condition \cref{th-CondNec} is efficient at proving that an initial state is not null-controllable: anything that is not smooth on $\T\setminus \omega$ is not null-controllable, but anything that is \emph{too smooth} (analytic) on the whole space is not null-controllable neither (\cref{th-analy-non-nc}).

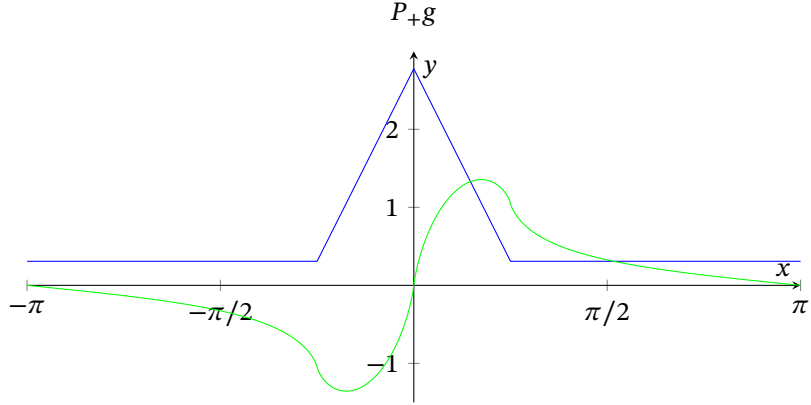
\begin{figure}
    \centering
    \begin{tikzpicture}
        \begin{axis}[
            title=$P_+g$,
            xlabel={$x$},
            ylabel={$y$},
            ymin=-1.5,
            ymax=3,
            axis lines=middle,
            width = 0.8\textwidth,
            height=0.3\textheight,
            xtick={-3.1415,-1.5708,0,1.5708,3.1415},
            xticklabels = {$-\pi$,$-\pi/2$,0,$\pi/2$,$\pi$},
            ytick={-1,0,1,2}
            ]
            \addplot [blue] table [x=x, y=re] {dilog.dat};
            \addplot [green] table [x=x, y=im] {dilog.dat};
        \end{axis}
    \end{tikzpicture}
    \caption{A plot of $\theta\mapsto P_+g(\eu^{\iu\theta})$ for $\theta_0 = \pi/4$. The real part is in blue, the imaginary part in green.}
    \label{fig-Pg}
\end{figure}

Giving a null-controllable initial state with a closed form expression is harder. Consider the case  $\omega = \{|z|=1,\ \lvert\arg(z)\rvert < a\}$. For $0<\theta_0\leq a$, consider $g$ a ``triangle'' function
\[
    g(\eu^{\iu \theta}) =
    \left\{
    \begin{IEEEeqnarraybox}[][c]{l'l}
         0,& -\pi<\theta < -\theta_0\\
         \theta_0 + \theta, & -\theta_0<\theta<0\\
         \theta_0-\theta,& 0<\theta<\theta_0\\
         0,& \theta_0<\theta<\pi
    \end{IEEEeqnarraybox}
    \right.
\]
Then straightforward (but tedious) computations show that
\[
 P_+g(z) = \theta_0^2 -2\Li(z) + \Li(\eu^{\iu\theta_0}z) + \Li(\eu^{-\iu \theta_0}z),
\]
where $\Li$ is the dilogarithm defined by $\Li(z) = \sum_{n>0} \frac{z^n}{n^2}$. We plot $P_+g$ in \cref{fig-Pg}. Since $g\in W^{1,2}_0(\omega)$, $P_+g\in \NCH(\omega)$. This is the simplest null-controllable initial state with a closed-form expression we could find.

A difficult question is to decide in general whether a given function $f_0\in H^2$ is the Riesz projection of a function $g\in L^2(\T)$ with $\supp(g) \subset \cl\omega$, and moreover extending holomorphically to $\hat{\C}\setminus \cl\omega$. Still we can make the following observation. Since $f_0$ extends holomorphically to $\Cinf\setminus \cl\omega$ (this can also be seen from \eqref{CauchyTrans}), we can deduce that the singular factor of $f_0$ has its boundary spectrum in $\omega$ (see \cite[Section A3]{NikVol1} for a general reference on factorization in Hardy spaces and \cite[Lemma A4.3.4]{NikVol1} for analytic continuation). %Additionally, since the singular inner factor and the outer factor are holomorphic outside $\T$, the Blaschke factor cannot have any zeros since otherwise $f_0$ would have a pole in the symmetric points. As a result $f_0=IF$ where $I$ is a singular inner function having its defining measure supported on $\omega$ and $F$ is an outer function which extends holomorphically through $\T\setminus\cl\omega$. Notice that $F=\exp(\frac{1}{2\pi}\int_{\T}\frac{e^{ix}+z}{e^{ix}-z}\log|h(e^{ix})|dx)$ is essentially the exponential of a Cauchy transform. 

We also mention Tumarkin's theorem~\cite[Theorem 5.3.1]{CMR06}, which we actually use to deduce \cref{th-CondNec} from \cref{th-NCH-holom-necessary}. So probably Tumarkin's theorem would not tell us more than \cref{th-NCH-holom-necessary}. We refer to \cite[Chapter 5]{CMR06} for related results. %, but it is unclear how they would be useful for our problem (to us, at least). However, here is a simple thing we can say:

Still, it is clear that if $f_0$ is the Riesz projection of $g\in L^2(\T)$ and $\supp(g)\subset\cl\omega$, then $g$ is unique. Indeed, if $P_+g=0$, then $P_-g=0$ on $\T\setminus \cl\omega$ which by standard uniqueness results in $H^2$~\cite[Theorem 17.18]{Rudin86} (or $\overline{H^2}$) implies that $g=0$. 

In the language of inverse problem, the problem ``recover $g\in L^2(\omega)$ from $P_+g$'' is identifiable. However, we do not have an inequality $\|g\|_{L^2(\omega)}\lesssim \|P_+g\|_{L^2}$: consider $\chi \in C_c^\infty(\omega)$ and $\chi_k(x) = \eu^{-\iu kx} \chi(x)$. Then $\widehat{\chi_k}(n) = \widehat \chi(n+k)$, hence $\|P_+ \chi_k\|_{L^2}^2 = 2\pi\sum_{n\geq k} |\widehat \chi(n)|^2\to 0$ as $k\to+\infty$ while $\|\chi_k\|_{L^2(\omega)}=c>0$. In the language of inverse problems, the problem above is not stable.

Let us finally discuss how the spaces (for the function $g$) that appear in \cref{th-CNS-W12,th-CondNec,th-CondSuff,th-CNS} compare to each other. It will be convenient to have a notation for the space of functions with different type of regularity/integrability for positive and negative frequencies.
\begin{definition}\label{def-XHD}
We assume that $\omega$ is a strict interval of $\T$ with endpoints $\zeta_1$ and $\zeta_2$. We set
     \begin{align*}
         X^{s+}_{t-}(\omega) & \coloneqq \{ g \in L^2(\T),\ \supp(g)\subset \cl\omega,\ P_+g \in W^{s,2}(\T),\ (I-P_+)g\in W^{t,2}(\T) \}\\
         \prescript{0}{}{X}^{s+}_{t-}(\omega) & \coloneqq X^{s+}_{t-}(\omega) \cap \left\{\int_{\T}\frac{|g(\eu^{\iu t})|^2}{|(\eu^{\iu t}-\zeta_1)(\eu^{\iu t}-\zeta_2)|}\diff t<\infty\right\}\\
         \widetilde X &\coloneqq \{g\in L^2(\omega),\ \dz P^{\Omega_T} g \in L^2(\Omega_T)\},
     \end{align*}
\end{definition}
Where $P^{\Omega_T}\colon L^2(\partial \Omega_T) \to h(\Omega_T)$ is the Poisson integral operator for $\Omega_T$. If $g\in L^2(\omega)$ and $f_0 = P_+g$, with these notations, the previous theorems can be rephrased as
\begin{IEEEeqnarray*}{+l+s}
    f_0 \in\NCH(\omega) \cap W^{1/2,2}(\T) \eq g \in X_{1/2-}^{1/2+}(\omega)& (\cref{th-CNS-W12})\\
    f_0 \in \NCH(\omega) \implies g\in X_{1/2-}^{0+}(\omega) & (\cref{th-CondNec})\\
    f_0\in \NCH(\omega) \eq g\in \widetilde X & (\cref{th-CNS})\\
    g \in \prescript{0}{}{X}_{1/2-}^{0+}(\omega) \implies f_0 \in \NCH(\omega) & (\cref{th-CondSuff})
\end{IEEEeqnarray*}

Thus, we have
\[
 X^{1/2+}_{1/2-}(\omega) \subseteq \prescript{0}{}{X}_{1/2-}^{0+}(\omega) \subseteq \widetilde X \subseteq {X}_{1/2-}^{0+}(\omega).
\]
See, e.g., \cite[Theorem 1.87]{Leoni23} for the first inclusion.
The space $X_{1/2-}^{1/2+}(\omega) = W^{1/2,2}_{00}(\omega)$ that appears in \cref{th-CNS-W12} is the simplest, or maybe the least unusual of the bunch. The triangle function $g$ above belongs in this space (and even in $W^{1,2}_0(\omega)$).

We can construct a function in $X^{0+}_{1/2-}(\omega) \setminus X^{1/2+}_{1/2-}(\omega)$ the following way: first consider $h(z) = -\sum_{n>0} z^n/n = \log(1-z)$. Then $P_+ h\notin W^{1.2,2}(\T)$ and $(I-P_+)h = 0 \in W^{1/2,2}(\T)$. Consider also $\chi \in C_c^\infty(\omega)$ such that $\chi =1$ on a neighborhood of $1$. Finally, set $g = \chi h$. Then, we have $g\in L^2(\omega)$, and writing $\widehat g(n) = \sum_{k<n} \widehat \chi(k)/(n-k)$, we see that $\widehat g(n)$ decays faster than any polynomial as $n\to -\infty$. In particular, $(I-P_+)g\in W^{1/2,2}(\T)$. On the other hand, $g\notin W^{1/2,2}(\T)$, as that would imply $h = (1-\chi)h + g \in W^{1/2,2}(\T)$, a contradiction. Thus $g$ is a function that satisfies the necessary condition of \cref{th-CondNec} that is not in $W^{1/2,2}(\T)$. Since the support of $g$ is compact in $\omega$, $g$ also satisfies the sufficient condition of \cref{th-CondSuff}.

Thus, in the inclusions chain above, we actually have $X^{1/2+}_{1/2-}(\omega) \subsetneq \prescript{0}{}{X}^{0+}_{1/2-}(\omega)$. We could not prove whether the other inclusions are strict or not.

\subsection{Bibliographical comments}
\paragraph{Control of the fractional heat equation}
The controllability of the fractional heat equation
\[
    (\partial_t+(-\Delta)^{\alpha}u)f=Bu
\]
has been studied since 2006 by Miller~\cite{Miller06} and Micu and Zuazua~\cite{MZ06}.

More precisely, Miller~\cite{Miller06} shows that the Lebeau-Robbiano control strategy (used to prove the null-controllability of the heat equation~\cite{LR95}) is also adapted to prove the null-controllability of the fractional heat equation when $\alpha>1/2$.

On their end, Micu and Zuazua use the moment method (to prove the null-controllability of the heat equation in dimension one~\cite{FR71}). They also remark that this strategy works for the fractional heat equation when $\alpha >1/2$.

In the case $\alpha \leq 1/2$, Miller, Micu and Zuazua (see also~\cite[Appendix]{DM12}) consider \emph{shaped controls}, i.e., the right-hand side of the equation is $u(t)h(x)$ (instead of $\mathds 1_\omega u(t,x)$ as we are considering in this article). Using M\"untz' theorem, they show that the fractional heat equation with $\alpha\leq 1/2$ is not null-controllable: for every control shape $h(x)$ and $T>0$ there exists an initial state that cannot be steered to $0$ in time $T$.

We emphasize that the kind of control they consider (shaped control) is different to the controls we are considering (internal control). On the surface, both kinds controls seem to give similar results in dimension one: if $\alpha >1/2$ the fractional heat equation is null-controllable with both, and if $\alpha \leq 1/2$, the fractional heat equation is not null-controllable, whatever kind of control (see \cites{Koenig20,AK24,Lissy22} and \cite[Theorem 4]{Koenig17} for lack of null-controllability with internal controls).

However, differences in controllability properties appear when we look more closely: for certain control shapes $h$, the space of null-controllable states might be $\{0\}$~\cite[Theorem 3.1]{MZ06}. On the other end, with internal controls, our results show that the space of null-controllable states is somewhat large (\cref{th-CondSuff,th-CNS-W12,th-dense}). Let us also mention that \Cref{th-analy-non-nc} was previously proved by the second author with a different proof~\cite[Théorème 2.5]{Koenig19}.

These differences of controllability properties indicate that our methods do not adapt for shaped control. Indeed, an important point of our analysis consists in identifying the right-hand side of the observability inequality and a Bergman space norm (\cref{th-MainEstim}), but nothing of this sort is possible for shaped controls (see \cite[eq.~(5.1)]{MZ06} for the observability inequality associated to shaped controls).

The case $\alpha =1/2$ (i.e., the half-heat equation) is interesting for another reason: in a sense, it can be embedded in the \emph{Baouendi-Grushin heat equation} (where $\omega\subset \Omega\coloneqq(-1,1)\times \T$),
\begin{equation*}%\tag{$E_{L^2}$}\label{eq-half-L2}
    \left\{\begin{IEEEeqnarraybox}[][c]{l"l}
        (\partial_t -\partial_y^2-y^2\partial_x^2)f(t,y,x) =  \mathbf{1}_\omega u(t,y,x),
        & (t,y,x) \in \R_+ \times (-1,1)\times T\\
        f(0,y,x) = f_0(y,x),& (y,x)\in (-1,1)\times \T,
    \end{IEEEeqnarraybox}\right.
\end{equation*}
This time it is controlled from the interior of the boundary of the domain (we have interchanged $x$ and $y$ here in order to keep $x\in\T$ to be consistant with the notation of the present article). Consider now $(v_{n,k})$, the eigenfunctions of $-\partial_y^2+(ny)^2$
(with Dirichlet boundary conditions), then $\Phi_{n,k}(y,x)=v_{n,k}(y)e^{inx}$ are eigenfunctions of $-\partial_y^2-y^2\partial_x^2$ and free solution can be expressed as $f(t,y,x)=\sum_{n,k}\eu^{-t\lambda_{n,k}} c_{n,k} \Phi_{n,k}(y,x)$. Since $\lambda_{n,0} = n + o(n)$~\cite[Theorem 4.23]{DS99}, it is natural to introduce the corresponding differential operator on $\T$ by
\[
    |D_x|(\sum_n a_ne^{inx})=\sum_na_n|n|\eu^{\iu nx}.
\]
With this in mind, after having established that the so-called half-heat equation $(\partial_t +|D_x|)f(t,x) = \mathds 1_\omega u(t,x)$ ($f$ being in a suitable subspace of $H^2$) is not null-controllable on $\omega$ in any time, the second author deduces that the Grushin equation is not null controllable in any time~\cite[Theorem 4]{Koenig17}.

This is the starting point of this paper. Indeed, knowing that the half-heat equation is not null-controllable in any time raises the question which states can be steered to zero.

\paragraph{Complex analysis tools for control theory}
Complex analysis tools have proven to be powerful in the study of several control problems, and our article is in the line of this tradition.

In particular, holomorphic function spaces and analytic function theory have shown to be efficient in describing reachable states in certain control problems. Let us mention the classical heat equation. For an interval $I\subset \R$, this is given by
\begin{equation*}
    %\label{HE}
    %\tag{HE}
    \left\lbrace
    \begin{aligned}
         & \frac{\partial y}{\partial t}(t,x)-		\frac{\partial^2 y}{\partial x^2} = 0  \qquad & t >0, \ x\in I, & \\
        %		&y(t,0)=u_0(t),\  \ y(t,\pi)=u_\pi(t) \qquad & t >0, &\\
         & y(0,x)= f(x) \qquad                                                                & x \in I,
    \end{aligned}
    \right.
\end{equation*}
Aikawa, Hayashi and Saitoh \cite{AHS90} were considering this problem when $I=\R_+$ and the equation is controlled in $x_0=0$ by an $L^2$-function in some finite time. As it turns out, the solution $y(T,\cdot)$ ($T>0$ fixed) can be expressed as a certain Laplace-type transform for which the range is characterized as a weighted Bergman space on the sector $\{z=x+iy\in\C:x>0,|y|<x\}$. %We recall that for a domain $\Omega\subset \C$, the Bergman space $A^2(\Omega)$ is the set of functions $f$ holomorphic on $\Omega$ such that $\|f\|_{A^2(\Omega)}^2=\int_{\Omega}|f|^2\diff A<+\infty$, where $A$ is planar Lebesgue measure (weighted versions are defined integrating against the corresponding weight). 

Another prominent situation is when $I$ is a finite interval, say $I=(0,\pi)$, and the heat equation is controlled from $x_0=0$ and $x_1=\pi$:
\begin{equation*}
    %\label{HE}
    %\tag{HE}
    \left\lbrace
    \begin{aligned}
         & \frac{\partial y}{\partial t}(t,x)-		\frac{\partial^2 y}{\partial x^2} = 0  \qquad & t >0, \ x\in (0,\pi), & \\
         & y(t,0)=u_0(t),\  \ y(t,\pi)=u_\pi(t) \qquad                                        & t >0,                 & \\
         & y(0,x)= f(x) \qquad                                                                & x \in (0, \pi),
    \end{aligned}
    \right.
\end{equation*}
where $u_0$ and $u_{\pi}$ are $L^2$-functions (Dirichlet control; other types of boundary conditions can be considered). It is a well known fact that this equation is null-controllable in any time, and thus the reachable states do not depend on time. Consider $u_0$ and $u_{\pi}$ on some finite interval $(0,T)$. In this situation, pioneering work by Fattorini and Russell \cite{FR71} showed that the reachable states extend to holomorphic functions on the square $D$ one diagonal of which is the interval $(0,\pi)$. Subsequent work focused on a more precise description of the reachable space. Martin, Rosier and Rouchon \cite{MRR16} proved that functions holomorphic in a disk with a certain radius are reachable. This was significantly improved by Dard\'e and Ervedoza \cite{DE18} who establish that functions holomorphic in a neighborhood of $D$ are reachable. Tuscnak, Hartmann and Kellay \cite{HKT20} then showed that the reachable space is sandwiched between two spaces of holomorphic spaces on $D$, the Smirnov space, which is a kind of companion space to the Hardy space, and the Bergman space. The final step was reached in \cite{HO21b} where the reachable space was characterized as the Bergman space on $D$. A central tool of this result was the so-called separation of singularities problem which is a purely complex analysis problem (see below), and which will show to be of relevance also in the problem we will consider here.

\paragraph{Null-controllable states for other equations}
Most paper on null-controllability either
\begin{itemize}
    \item prove that a control system is null-controllable, e.g., the heat equation~\cite{LR95,FI96}, Baouendi-Grushin equation in some context~\cite{BCG14}, fractional heat with exponent $>1/2$~\cite{Miller06,Miller10,DM12,AM22b} and many more,
    \item or prove it is not null-controllable, e.g., fractional heat equation with exponent $\leq 1/2$ \cite{Koenig17,Koenig20,AK24} Baouendi-Grushin in other contexts~\cite{BCG14,Koenig17} and many more.
\end{itemize}
In the second case, one could ask which initial conditions can be steered to $0$, but this is less studied. We are aware of the following paper that tackle this question:
\begin{itemize}
    \item The heat equation (possibly with some quadratic potential) on the half-line \cite{MZ01,DE19}, where no non-zero initial state can be steered to $0$,
    \item The Baouendi-Grushin equation controlled on a symmetric vertical band \cite[Theorem 1.4]{BMM15}, where initial states that are regular enough can be steered to $0$ in a smaller time than other initial states,
    \item the heat equation on $\R$ with ``less and less control at $\infty$'' \cite{CMV04}, where initial states in a suitably weighted $L^2$ space can be steered to $0$.
\end{itemize}
Our article is an addition to this short list.

\subsection{Notations}
We will introduce notations along the way, but to help the reader, we collect here the recurring notations that appear in this article.
\begin{IEEEeqnarray*}{l"s}
    \IEEEeqnarraymulticol{2}{s}{\emph{General notations}}\\
    A\lesssim B& $\exists C>0,\ A\leq CB$ ($A$, $B$ depend on some parameters)\\
    A\asymp B& $A\lesssim B$ and $B\lesssim A$\\
    \langle f, g\rangle_H& Scalar product in the Hilbert space $H$, left-linear and right-antilinear\\
    \IEEEeqnarraymulticol{2}{s}{\emph{Sets}}\\
    \Cinf & Riemann sphere ($\C\cup \{\infty\}$)\\
    \D & Unit disk $\{z\in \C,\ |z|<1\}$\\
    \T & Complex unit circle $\{z\in \C,\ |z|=1\}$\\
    \Dext & Extended exterior of the unit disk $\{z\in \C,\ |z|>1\}\cup \{\infty\}$\\
    \Omega & domain in $\C$\\
    \overline E & When $E \subset \C$, $\{\overline z,\ z\in E\}$\\
    \cl E & Closure of $E$\\
    \distance(E,F)& Distance between two subsets of $\C$\\
    \IEEEeqnarraymulticol{2}{s}{\emph{Spaces}}\\
    L^2(\T) & Lebesgue space on $\T$\\
    L^2(\Omega) & Lebesgue space on $\Omega$ with respect to area measure\\
    h(\Omega) & Space of harmonic functions on $\Omega$\\
    \Hol(\Omega) & Space of holomorphic functions on $\Omega$\\
    H^2 & Hardy space ($H^2 = H^2(\D) = P_+ L^2(\T)$, see below for $P_+$)\\
    A^2(\Omega) & Bergman space (holomorphic functions in $L^2(\Omega))$\\
    \Dirichlet(\Omega)& Dirichlet space (holomorphic functions with derivative in $L^2(\Omega))$\\
    W^{k,p},\ W^{k,p}_0,\ W^{k,p}_{00} & Sobolev spaces\\
    \NCH(\omega,T),\ \NCL(\omega,T) & Spaces of null-controllable initial states (\cref{def-NC})\\
    \IEEEeqnarraymulticol{2}{s}{\emph{Operators}}\\
    P_+& Riesz projector $L^2(\T)\to L^2(\T)$ (projection
    on positive frequencies)\\
    P_-& $I-P_+:L^2(\T)\to L^2(\T)$ ($I$ identity operator)\\
    S(t)& Half-heat semigroup (\cref{eq-semigroup})\\
    P^\Omega &Poisson kernel (and associated integral operator $L^2(\partial\Omega) \to h(\Omega)$)\\
    P^\eu& Poisson Kernel on $\Dext$\\
    \dz, \dzb& Wirtinger derivatives $(\partial_x-\iu \partial_y)/2$ and $(\partial_x+\iu\partial_y)/2$ respectively\\
    \FT & Input-to-output map (\cref{def-FT})\\
    \IEEEeqnarraymulticol{2}{s}{\emph{Measures}}\\
    \lvert \diff z|& Linear ($1$-dimensional Hausdorff) measure \\
    \diff A& Area ($2$-dimensional Lebesgue) measure\\
\end{IEEEeqnarray*}

\section{The \texorpdfstring{$H^2$}{H²} problem}
\subsection{The observability inequality}\label{sec-obs}
Let us start by defining the \emph{input-to-output map}; a usual tool in control theory.

\begin{definition}\label{def-FT}
    Let $T>0$. We define the input-to-output map associated to the $L^2$ control system~\eqref{eq-half-L2} $\FT\colon L^2(0,T;\ L^2(\omega)) \to L^2(\T)$ by
    \[
    \forall u \in L^2(0,T;\ L^2(\omega)),\ 
    \FT u(\eu^{\iu x})
    \coloneqq\int_0^{T}S(T-s)(\mathds 1_{\omega}u(s,\cdot))(\eu^{\iu x})\diff s.
    \]
    
    We also define the input-to-output map associated to the $H^2$ control system~\eqref{eq-half-H2} $\FTp\colon L^2(0,T;\ L^2(\omega)) \to H^2$ by
    \[
    \forall u \in L^2(0,T;\ L^2(\omega)),\ 
    \FTp u(\eu^{\iu x})
    \coloneqq\int_0^{T}S(T-s)(P_+\mathds 1_{\omega}u(s,\cdot))(\eu^{\iu x})\diff s.
    \]
\end{definition}
These definitions correspond to the integral term that appears in Duhamel's formula for the control systems. Notice also that $\FTp = P_+\FT$.

We now prove \cref{th-MainEstim} by combining the change of variables $z = \eu^{-t+\iu x}$ with the equivalence between null-controllability and observability.
\begin{delayedproof}{th-MainEstim}
    By definition of $S(T)$, the solution $f$ of the $H^2$-control system~\eqref{eq-half-H2} is
    \[
        f(T,\cdot) = S(T) f_0 + \FTp u.
    \]
    Hence, a function $f_0$ is null controllable at time $ T >0$ if and only if $S(T)f_0\in \Rg \FTp$. Consider the operators $C_2\colon f_1 \in \C f_0\mapsto S(T) f_1\in H^2$ and $C_3 = \FTp \colon L^2([0,T]\times \omega) \to H^2$. Then, $f_0$ is null-controllable at time $T$ if and only if $\range(C_2) \subset \range(C_3)$.

    According to a standard duality lemma~\cite[Lemma~2.48]{Coron07}, this is equivalent to the existence of $C_{T,f_0}'>0$ such that
    \begin{equation}\label{eq-obs-h2-1}
        \forall g\in H^2,\ \|C_2^* g\|_{H^2}\leq C_{T,f_0}'\|C_3^*g\|_{L^2(0,T ,L^2(\omega))} =  C_{T,f_0}'\|(\FTp)^*g\|_{L^2(0,T ,L^2(\omega))},
    \end{equation}
    and if this holds, there exists a continuous linear map $C_1\colon \C f_0 \to L^2([0,T]\times\omega)$ such that $C_2 = C_3 C_1$ and $\|C_1\| \leq C_{T,f_0}'$. This relation reads $S(T) f_0 = \FTp C_1 f_0$, that is, $-C_1 f_0$ steers $f_0$ to $0$ in time $T$. Thus, if the estimate~\eqref{eq-obs-h2-1} holds, we can choose the control $u$ that steers $f_0$ to $0$ such that $\|u\|_{L^2([0,T]\times\omega)} \leq C_{T,f_0}' \|f_0\|_{H^2}$.

    To compute the left and right-hand sides, we first see from the definition of the semi-group $S$ (see \eqref{eq-semigroup}) that if $g\in H^2(\D)$ (identified with $H^2(\T)$), then $S(t)^*g(\eu^{\iu x}) = S(t)g(\eu^{\iu x}) = g(\eu^{\iu x-t})$.

    Let us compute the left-hand side. Set $\iota_{f_0}$ the injection map from $\C f_0$ to $H^2$. Then, $C_2 = S(T)\circ \iota_{f_0}$. Routine computations show that $\iota_{f_0}^*$ is the orthogonal projection onto $\C f_0$, that is
    \[
        \iota_{f_0}^* g = \|f_0\|_{H^2}^{-2}\langle g,f_0\rangle_{H^2} f_0.
    \]
    Thus,
    \begin{equation}\label{eq-obs-lhs}
        C_2^* g = \iota_{f_0}^*\circ S(T)^* g 
        = \iota_{f_0}^*\circ S(T) g 
        = \iota_{f_0}^*( g(\eu^{-T}\cdot) )
        = \|f_0\|_{H^2}^{-2} \langle g(\eu^{-T}\cdot),f_0\rangle_{H^2} f_0.
    \end{equation}
    Thus the observability inequality~\eqref{eq-obs-h2-1} is equivalent to
    \begin{equation}\label{eq-obs-h2-2}
        \forall g\in H^2,\ |\langle g(\eu^{-T}\cdot),f_0\rangle_{H^2}|\leq C_{T,f_0}'\|f_0\|_{H^2}\|(\FTp)^*g\|_{L^2(0,T ,L^2(\omega))}.
    \end{equation}

    Let us now compute the right hand side of \cref{eq-obs-h2-1}. First, routine computations~\cite[Lemma~2.47]{Coron07} show that
    \begin{equation}\label{adjoint}
        (\FTp)^*g(t) = \mathds 1_\omega^* P_+S(T -t)^*g,
    \end{equation}
    where $\mathds 1_\omega^*\colon L^2(\T) \to L^2(\omega)$ is the restriction operator.
    Thus, the right-hand side of the observability inequality~\eqref{eq-obs-h2-1} is,
    \begin{IEEEeqnarray*}{rCl}
        \|(\FTp)^*g\|_{L^2(0,T ;L^2(\omega))}^2
        &=&\int_0^{T } \|\mathcal(\FTp)^*g(t)\|_{L^2(\omega)}^2\diff t
        =\int_0^{T }\int_{\omega}|g(\eu^{\iu x -(T-t)})|^2 \diff x\diff t\\
        &=&\int_0^{T }\int_{\omega}|g(\eu^{\iu x -t})|^2 \diff x\diff t.
        \intertext{We make the change of variables $z=\eu^{\iu x-t}$, that satisfies $\diff A(z) = |z|^{-2} \diff t\diff x$ ($\diff A$ is the area measure on $\C \simeq \R^2$). The image of $[0,T ]\times \omega$ by this change of variables is $\Omega_{T}^*=\{z\in \C \colon \arg(z)\in\omega, \eu^{-T }<|z|<1\}$ (see \cref{fig-Omega}). Hence,}
        \|(\FTp)^*g\|_{L^2(0,T ;L^2(\omega))}^2
        &=& \int_{\Omega_T^* } |g(z)|^2|z|^{-2} \diff A(z).\IEEEyesnumber\label{eq-obs-rhs}
    \end{IEEEeqnarray*}
    Since $1\leq|z|^{-2} \leq \eu^{2T}$, plugging this in the right-hand side of~\eqref{eq-obs-h2-2} proves that the observability inequality~\eqref{eq-obs-h2-1} (hence the null-controllability of $f_0$) is equivalent to
    \begin{equation}\label{eq-obs-h2}
        |\langle g(\eu^{-T}\cdot),f_0\rangle_{H^2}|
        \leq C_{T,f_0}'' \|g\|_{A^2(\Omega_T^* )}
    \end{equation}
    for every $g\in H^2$. 
    %In addition, if the inequality~\eqref{eq-obs-h2-1} holds, we can choose $C_{T,f_0}' \leq \eu^{T}\|f_0\|^2_{H^2} C_{T,f_0}$.
    In addition, if the inequality~\eqref{eq-obs-h2} holds, then, the inequality~\eqref{eq-obs-h2-2} (and hence \eqref{eq-obs-h2-1}) holds with $C_{T,f_0}' \leq \|f_0\|_{H^2}^{-1} C_{T,f_0}''$. In this case, we can choose the control $u$ that steers $f_0$ to $0$ such that $\|u\|_{L^2([0,T]\times\omega)}\leq C_{T,f_0}' \|f_0\|_{H^2} \leq C_{T,f_0}''$.
    
    Since polynomials are dense in $H^2$~\cite[Theorem 3.3]{DS04}, and since $H^2 \subset A^2(\D) \subset A^2(\Omega_T^*)$ with continuous inclusions, we can replace ``for every $g\in H^2$'' above by ``for every $g\in\C[X]$''.
    
    Finally, the change of function $\widetilde g(z) = g(\eu^{-T }z)$ proves that the inequality~\eqref{eq-obs-h2} can be written as
    \begin{equation}
        \forall g\in \C[X],\ |\langle f_0,g\rangle_{H^2}|
        \leq \eu^{-T}C_{T,f_0}'' \|g\|_{A^2(\Omega_{T })},
    \end{equation}
    where $\Omega_{T}=\eu^T \Omega_T^* $ is indeed the one defined in the statement of \Cref{th-MainEstim}.
\end{delayedproof}

\subsection{Independence of null-controllable states with respect to time}\label{sec-T}

We first prove that $\NCH(\omega,T)$ does not depend on $T$. In fact, we are able to estimate the norm of the control as $T\to 0$:
\begin{theorem}\label{th-NCH-T-2}
    Let $T>0$ and $f_0\in\NCH(\omega,T)$. There exists $C>0$ such that for every $T'>0$, there exists $u\in L^2([0,T']\times\omega)$ with
    \begin{itemize}
        \item $S(T')f_0 + \FTp[T'] u = 0$,
        \item $\|u\|_{L^2([0,T']\times\omega)} \leq C\max(1,T'^{-3/2})$.
    \end{itemize}
\end{theorem}

The first item says that the solution $f$ of the $H^2$ control system~\eqref{eq-half-H2} satisfies $f(T',\cdot) = 0$, i.e., $f_0\in \NCH(\omega,T')$. Hence, \cref{th-NCH-T-2} implies \cref{th-NCH-T}.

A central ingredient will be a result on separation of singularities:

\begin{theorem}[Corollary 1.2 \cite{HO21b}]\label{thm:SepSing}
    Let $1<p<+\infty$. Let $\Omega_1$ and $\Omega_2$ be open sets of $\C$ such that $\emptyset\neq \Omega_1\cap \Omega_2$ is bounded. If $\distance(\Omega_1\setminus\Omega_2,\Omega_2\setminus\Omega_1)>0$, then $A^p(\Omega_1\cap\Omega_2)=A^p(\Omega_1)+A^p(\Omega_2)$. %Moreover, there exist bounded linear operators $B_k:A^p(\Omega_1\cap\Omega_2)\to A^p(\Omega_k)$, $k=1,2$, such that for every $f\in A^p(\Omega_1\cap\Omega_2)$ we have $f=B_1f+B_2f$.
\end{theorem}

We actually need a little bit more information than this reference. Namely, we want the decomposition $A^2(\Omega_1\cap \Omega_2) \ni g = g_1+g_2 \in A^2(\Omega_1) + A^2(\Omega_2)$ to be done in a linear continuous way, with an upper bound on the norm of the resulting operator. All of that is done by looking carefully at the proof of the previous theorem, which we will do now:
\begin{theorem}\label{th-SepSing-est}
    Let $R>0$. Let $\Omega_1,\Omega_2\subset B(0,R)$ open such that $\distance(\Omega_1\setminus\Omega_2, \ \Omega_2\setminus\Omega_1)>0$. There exists $B_i\colon A^p(\Omega_1\cap\Omega_2) \to A^p(\Omega_i)$ ($i=1,2$) linear continuous such that
    \begin{itemize}
        \item $\forall g \in A^p(\Omega_1\cap \Omega_2)$, $g = B_1g + B_2g$,
        \item there exists $C_R>0$ depending only on $R$ such that for all $g\in A^p(\Omega_1\cap \Omega_2)$, 
        \[
            \|B_i g\|_{A^p(\Omega_i)}\leq \frac{C_R}{\distance(\Omega_1\setminus \Omega_2,\ \Omega_2\setminus \Omega_1)} \|g\|_{A^p(\Omega_1\cap\Omega_2)}.
        \]
    \end{itemize}
\end{theorem}
\begin{proof} The idea of the proof is essentially the same as in \cite{HO21b}. The starting point is the construction of a smooth partition of unity $\{1-\chi,\chi\}$ subordinate to $\{\Omega_1,\Omega_2\}$:  $\chi \in C^\infty(\C)$ such that $\chi = 0$ on a neighborhood of $\cl{(\Omega_1\setminus \Omega_2)}$ and $\chi=1$ on a neighborhood of $\cl{(\Omega_2 \setminus \Omega_1)}$ with $0\leq \chi \leq 1$. The key observation, which is the novelty here, is that $\chi$ additionally satisfies
\begin{equation}\label{eq-cutoff-est}
    \|\nabla \chi\|_{L^\infty(\Omega_1\cup\Omega_2)} \leq \dfrac{C}{\distance(\Omega_1\setminus \Omega_2,\ \Omega_2\setminus \Omega_1)}.
\end{equation}
The construction is standard, but we detail it in \cref{sec-cutoff} for the reader's convenience.

Now, given $g\in A^p(\Omega_1\cap \Omega_2)$, set $h = g\dzb \chi$ (initially defined on $\Omega_1\cap \Omega_2$ and extended by $0$) and $u = \frac1z\ast h$, the operators $B_i$ are defined as follows:
\begin{equation}
    \begin{aligned}
        B_1 g &= \chi g -u\\
        B_2 g &= (1-\chi) g +u,
    \end{aligned}
\end{equation}
where the functions on the right hand side are extended trivially outside $\Omega_1\cap\Omega_2$ (these extensions are obviously smooth on $\Omega_1$ and $\Omega_2$ respectively).

%\step{Conclusion}
%We have
By definition $B_1g + B_2 g = g$. In addition, since $1/z$ is the fundamental solution for the $\dzb$ operator, we get that $\dzb u = g\dzb\chi$. Hence $\dzb B_1g = 0$ and $\dzb B_2 g = 0$. That is, $B_1g$ and $B_2g$ are holomorphic.

We finally obtain the estimates of $\|B_i\|$: %This is a consequence of the estimates of the first step ($0\leq\chi \leq1$ and $|\nabla \chi|\leq C/\epsilon$):
\begin{align*}
    \|B_1 g\|_{L^p(\Omega_1)}
    &\leq \|g\|_{L^p(\Omega_1 \cap \Omega_2)} + \big \|\frac1z \big\|_{L^1(B(0,2R))} \|g\dzb\chi\|_{L^p(\Omega_1\cap \Omega_2)}\\
    &\leq \left(1+\big \|\frac1z \big\|_{L^1(B(0,2R))}\times\dfrac{C}{\distance(\Omega_1\setminus \Omega_2,\ \Omega_2\setminus \Omega_1)}\right) \|g\|_{L^p(\Omega_1 \cap \Omega_2)},
\end{align*}
and similarly for $B_2$.
\end{proof}

We can now prove \cref{th-NCH-T-2}.
\begin{delayedproof}{th-NCH-T-2}
    Suppose %$E_{H^2})$ is stearable to zero at time $T>0$, and 
    that the function $f_0$ is steerable to zero in time $T>0$. According to \cref{th-MainEstim}, this means that there exists $C_{T,f_0}>0$ such that for every polynomial $g$,
    \begin{equation}\label{MainEst2}
        |\langle f_0,g\rangle_{H^2}| \leq C_{f_0,T}\|g\|_{A^2(\Omega_T)}.
    \end{equation}
    
    Now, if $T'>T$, then $\Omega_T\subset \Omega_{T'}$, so that trivially $\|g\|_{A^2(\Omega_T)}\leq \|g\|_{A^2(\Omega_{T'})}$ for $g\in A^2(\Omega_{T'})$, which in particular holds for polynomials. We thus get \eqref{MainEst2} for the norm in $A^2(\Omega_{T'})$ with the same constant $C_{f_0,T}$, showing that $f_0$ is null controllable in time $T'>T$ with control of norm at most $C_{f_0,T}$.

    Consider the situation when $0<T'<T$. 
    We have to check the observability inequality~\eqref{MainEstim} for every polynomial $g$ knowing that it holds for functions in $A^2(\Omega_T)$. Set $\widehat{\Omega}_{T'}=\eu^{T'/2}\D\cup \Omega_{T'}$, which compactly contains the closure of $\D$, see \cref{fig-Omega-hat}.

    \begin{figure}
        \begin{minipage}[c]{0.5\textwidth}
            \begin{center}
            \tikzmath{\q=1.2; \Rad=1.7; \Radp=1.4; \a=10; \b=60; \acomp = 360+\a; \abmid = (\a+\b)/2; \Rmidp = (1+\Radp)/2; \Rmid = (\Rad+\Radp)/2; \fr = 0.6/\Rad; \frp = 0.6/\Radp;}
\begin{tikzpicture}[scale=2, every node/.style={fill=white, fill opacity = 0.3, inner sep=1pt, text=black, text opacity=1}]
        
    \draw[fill=Omegac, fill opacity = 0.5] (\a:1) -- (\a:\Rad)
        arc[start angle=\a, end angle=\b, radius=\Rad]
        -- (\b:1)
        arc[start angle=\b, end angle=\a, radius=1];
        
    \draw[thick, diskc!50!black, fill=diskc, fill opacity=0.4] (\a:\q) -- (\a:\Radp)
        arc[start angle=\a, end angle=\b, radius = \Radp]
        -- (\b:\q)
        arc[start angle=\b, end angle=\acomp, radius = \q];
    
    \draw[|-|, ultra thick, omegac] (\a:1) arc[start angle=\a, end angle = \b, radius = 1] node[pos = 0.5, below left]{$\omega$};
    
    \node at (\abmid:\Rmid) {$\Omega_T$};
    \node at (\abmid:\Rmidp) {$\Omega_{T'}$};
    % \draw[->] (-35:\Radp) node[below]{$\Omega_{T'}$} -- (\abmid:\Rmidp);
    
    \node at (-135:0.5){$\widehat{\Omega}_{T'}$};
    \draw[gray, thin] (0,0) circle[radius=1];
    \node[above] at (-45:1) {$\T$};
    
    \draw[<->] (0,0) -- ($(\b+70:\q)$) node[pos=0.5, auto=right]{$\eu^{T'/2}$};
    \draw[<->] (0,0) -- ($(4*\a/5+\b/5:\Rad)$) node[pos=\fr, auto=right]{$\eu^{T}$};
    \draw[<->] (0,0) -- ($(\a/5+4*\b/5:\Radp)$) node[pos=\frp, auto=left]{$\eu^{T'}$};
    
    \draw[->] (-1.2,0) -- (1.9,0);
    \draw[->] (0,-1.2) -- (0,1.4);
\end{tikzpicture}
            \end{center}
        \end{minipage}\hfill%
        \begin{minipage}[c]{0.4\textwidth}
            \caption{Illustration of $\widehat{\Omega}_{T'}$ and $\Omega_T$ when $T'<T$. The partial ring $\Omega_T$ is in red (including the reddish-yellow part). $\widehat\Omega_{T'}$ is in yellow (also including the reddish-yellow part). The intersection of $\Omega_T$ and $\widehat \Omega_{T'}$ is $\Omega_{T'}$ colored in reddish-yellow.}
            \label{fig-Omega-hat}
        \end{minipage}
    \end{figure}
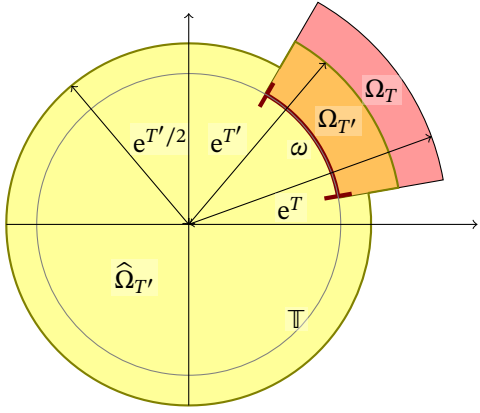

    Then
    \[
        \Omega_{T'}=\widehat{\Omega}_{T'}\cap \Omega_T.
    \]
    By construction, $\widehat{\Omega}_{T'}\setminus \Omega_T$ and $\Omega_T\setminus \widehat{\Omega}_{T'}$ are at a strictly positive distance $\eu^{T'}-\eu^{T'/2}\geq T'/2$, and so, using Theorem \ref{thm:SepSing}, we have
    \[
        A^2(\Omega_{T'})=A^2(\widehat{\Omega}_{T'})+A^2(\Omega_T).
    \]
    Moreover, for every $g\in A^2(\Omega_{T'})$, we have $g=B_1g+B_2g$, where $B_1g\in A^2(\widehat{\Omega}_{T'})$ and $B_2g\in A^2(\Omega_T)$ and for some $C_T>0$,
    \begin{align*}
    \|B_1 g\|_{A^2(\widehat{\Omega}_{T'})} &\leq \frac{C_T}{T'}\|g\|_{A^2(\Omega_{T'})},\\
    \|B_2 g\|_{A^2(\Omega_T)} &\leq \frac{C_T}{T'}\|g\|_{A^2(\Omega_{T'})}.
    \end{align*}
    
    Pick now a polynomial $p$, and let  $p=B_1p+B_2p$ be a decomposition of $p$. Note that $B_1p$ and $B_2p$ need not to be polynomials. Define an anti-linear functional by
    \[
        \ell_{f_0}(p)=\langle f_0,p\rangle_{H^2}, \quad p\text{ polynomial.}
    \]
    We have to show that $|\ell_{f_0}(p)|\lesssim T'^{-3/2}\|p\|_{A^2(\Omega_{{T'}})}$ for every polynomial $p$.
    
    Now, since $f_0$ is null controllable in time $T$, by the observability inequality \eqref{MainEstim} we have for every polynomial $p$
    \[
        |\ell_{f_0}(p)|=|\langle f_0,p\rangle_{H^2}|\leq C_{T,f_0}\|p\|_{A^2(\Omega_T)}.
    \]
    Since polynomials are dense in $A^2({\Omega}_{T})$ we can extend $\ell_{f_0}$ to a unique bounded functional $\ell_1$ defined on the whole space $A^2({\Omega}_{T})$ (note that all sets on which we consider Bergman spaces in this paper are so-called Carath\'eodory domains, see \cite[Definition 1.6, Section 18.1]{Conway95} which gives density of polynomials). We have $\|\ell_1\|\leq C_{T,f_0}$.
    
    On the other hand, since $\cl\D\subset D(0,\eu^{T'/2}) \subset \widehat{\Omega}_{T'}$, % we have $A^2(\widehat{\Omega}_{T'})\subset H^2$, and
    an elementary computation leads to
    \begin{align*}
    \|g\|_{H^2}^2
    &=\sum_{n\ge 0}|\widehat{g}(n)|^2\\
    &=\frac{2}{T'}\sum_{n\ge 0}|\widehat{g}(n)|^2\frac{T'}{2}\\
    &\le \frac{2}{T'}\sum_{n\ge0}|\widehat{g}(n)|^2\frac{\eu^{T'(n+1)}-1}{2n+2}\\
    &=\frac{1}{T'}\|g\|^2_{A^2(\eu^{T'/2}\D\setminus \cl\D)}\\
    &\le \frac{1}{T'}\|g\|_{A^2(\widehat{\Omega}_{T'})}^2
    %\|g\|_{H^2}^2 
    %&= \sum_{n\geq 0} |\widehat g(n)|^2\\
    %&= \sum_{n\geq 0} |\widehat g(n)|^2 \frac{\eu^{Tn}}{n+1} (n+1)\eu^{-Tn}\\
    %& \leq \left(\max_{\lambda>0} \lambda\eu^{-\lambda}\right) \frac{\eu^{T'}}{T'} \sum_{n\geq 0} |\widehat g(n)|^2 \frac{\eu^{-T' n}}{n+1}\\
    %&\lesssim  T'^{-1}\|g\|_{A^2(\widehat{\Omega}_{T'})}^2
    \end{align*}
    Hence, the Cauchy-Schwarz inequality yields
    \begin{align*}
        |\ell_{f_0}(p)|^2&=|\langle f_0,p\rangle_{H^2}|^2
        = \left|\int_{\T} f_0(\eu^{\iu x})\overline{p(\eu^{\iu x})} \diff x\right|^2\le \|f_0\|_{H^2}^2\int_{\T}|p(\eu^{\iu x})|^2\diff x\\
        &\lesssim T'^{-1}\|p\|_{A^2(\widehat{\Omega}_{T'})}^2.
    \end{align*}
    Analogously, since polynomials are dense in $A^2({\Omega}_{T})$ we can extend $\ell_{f_0}$ to a unique functional $\ell_2$ defined on the whole space $A^2(\widehat{\Omega}_{T'})$, and $\|\ell_2\|\lesssim T'^{-1/2}\|f_0\|_{H^2}$.
    
    Since $A^2(\Omega_T)\cap A^2(\widehat{\Omega}_{1})$ ($=A^2(\Omega_T\cup \widehat{\Omega}_{1})$) is dense in $A^2(\Omega_{T'})$ (they all contain densely polynomials),
    we easily check that $\ell_1 = \ell_2$ on $A^2(\Omega_T)\cap A^2(\widehat{\Omega}_{1})$. Hence the formula 
    \[
        \forall (g_1,g_2) \in A^2(\Omega_T)\times A^2(\widehat{\Omega}_{T'}),\ \ell(g_1 + g_2) \coloneqq \ell_1(g_1) + \ell_2( g_2)
    \]
    defines uniquely a functional on $A^2(\Omega_T)+A^2(\widehat{\Omega}_{1})$. Clearly, by uniqueness, $\ell=\ell_1$ on $A^2(\Omega_T)$ and $\ell=\ell_2$ on $A^2(\widehat{\Omega}_{T'})$. Hence, for every polynomial $p$
    \[
        |\ell_{f_0}(p)|=|\ell(p)|
        =|\ell(B_1p+B_2p)|\le
        |\ell(B_1p)|+ |\ell(B_2p)|=|\ell_1(B_1p)|+ |\ell_2(B_2p)|.
    \]

    Now, since the $\ell_i$ are bounded on the underlying spaces we get
    \[
        |\ell_{f_0}(p)|\lesssim \|B_1p\|_{A^2(\Omega_T)}+ T'^{-1/2}\|B_2p\|_{A^2(\widehat{\Omega}_{T'})}
        \lesssim T'^{-3/2}\|p\|_{A^2(\Omega_{{T'}})},
    \]
    which concludes the proof.
\end{delayedproof}

\subsection{Necessary conditions and sufficient conditions \emph{via} complex analysis}\label{sec-h2}
We next prove the extension result. Since we have established the independence of null-controlla\-bi\-lity in time, we will henceforth no longer use the index $T$ in $\Omega_T$ but simply write $\Omega$ to not overcharge notation. Recall that given $\Omega=\{r\eu^{\iu t},\  \eu^{\iu t}\in\omega,\ 1<r<\eu^{T}\}$ (domain outside $\D$) we denote $\Omega^*=\{r\eu^{\iu t},\ \eu^{\iu t}\in\omega,\ \eu^{-T}<r<1\}$ (domain inside $\D$).
\begin{delayedproof}{th-NCH-holom-necessary}
    \step{Holomorphic extension}
    For $u\in \D$, let $k_u(z)=\frac{1}{1-\overline{u}z}$ be the reproducing kernel of $H^2$, which means that for every $f\in H^2(\D)$, $\langle f,k_u\rangle_{H^2}=f(u)$. The  function $k_u$ extends to a holomorphic function in $\Omega$ for every $u\in \C\setminus \Omega^*$, and it is clear that for every $u\in \C\setminus \cl{(\Omega^*)}$, we have $k_u\in A^2(\Omega)$. Recall the anti-linear functional introduced previously %on $A^2(\Omega^*)$ by
    \[
        \ell_{f_0}(p)=\langle f_0,p\rangle_{H^2}, \quad p\text{ polynomial.}
    \]
    Polynomials being dense in $A^2(\Omega)$, %, see \cite{Bourdon87}. 
    $\ell_{f_0}$ is well defined and continuous on $A^2(\Omega)$.
    Since $k_u\in A^2(\Omega)$ for $u\in \C\setminus \cl{\Omega}$, we can define
    \[
        \phi(u)=\ell_{f_0}(k_u).
    \]

    Now, the function $u\in \C\setminus \cl{(\Omega^*)} \mapsto k_u \in A^2(\Omega)$ is obviously anti-analytic, and so
    the function $\phi$ is analytic on $\C\setminus \cl{(\Omega^*)}$.
    Since $k_u$ is the reproducing kernel of $H^2$, we get for $|u|<1$ and $u\notin \cl{(\Omega^*)}$, $\phi(u) = \ell_{f_0}(u) = \langle f_0,k_u\rangle_{H^2} = f_0(u)$. Hence, since $f_0$ is holomorphic in $\D$, $f_0$ extends holomorphically to $\C\setminus\omega^{cl}$.

    Another consequence of the continuity of $\ell_{f_0}$ on $A^2(\Omega)$ is the existence of $\psi\in A^2(\Omega)$ such that $\ell_{f_0}(p)=\langle\psi,p\rangle_{A^2(\Omega)}$ and hence
    \[
        f_0(u)=\ell_{f_0}(k_u)=\langle\psi,k_u\rangle_{A^2}.
    \]

    \step{Estimate for $|u|\to \infty$} Since $\|k_u\|_{A^2(\Omega)} \leq C/|u|$ as $|u| \to, \infty$, we get
    \[
        |f_0(u)|\leq C\|\ell_{f_0}\|\frac1{|u|}.
    \]
    This proves that %$f_0$ extends holomorphically to $\Cinf$ with 
    extends to $\infty$ with $f_0(\infty) = 0$.

    \step{Regularity}
    We denote the reproducing kernel of $A^2(\Omega)$ by $k_u^{\Omega}$~\cite[\S1.1]{DS04}. Recall the following usual transformation formula of the reproducing kernel (see for instance \cite[\S1.3,~Theorem~3]{DS04})
    \begin{equation}\label{ChangeVarKernel}
        k^{\Omega_1}_{\lambda}(z)=k^{\Omega_2}_{\varphi(\lambda)}(\varphi(z))\varphi'(z)\overline{\varphi'(\lambda)},
        \quad z,\lambda\in\Omega_1,
    \end{equation}
    where $\varphi\colon\Omega_1\to\Omega_2$ is the conformal mapping from $\Omega_1$ to $\Omega_2$. Denoting as before by $\Deinf=\{z\in\Cinf\colon|z|>1\}$
    the outer disk, we get with $\varphi\colon\Dext\to\D$, $\varphi(z)=1/z$,
    \[
        k^{\Deinf}_{u}(z)=\frac{1}{(1-\overline{u}z)^2},\quad
        |u|>1,|z|>1.
    \]
    (Observe that this is the same formula as for $k_u^{\D}$, $|u|<1$.)
    We denote $P_{A^2(\Deinf)}\colon h\in L^2(\Deinf) \to L^2(\Deinf)$ the orthogonal projection onto $A^2(\Deinf)$, which is given by $P_{A^2(\Deinf)}f(u) =  \langle f,k^{A^2}_u\rangle_{L^2(\Deinf)}$.

    Remark that $\dzb[u] k_u(z) = z k^{\Deinf}_u(z)$. Hence, if we differentiate the formula $f_0(u) = \langle \psi, k_u\rangle_{A^2(\Omega)}$ in $u$, $|u|>1$, we get by Lebesgue's dominated convergence theorem
    \begin{align*}
        f_0'(u)
         & = \langle \psi, \dzb[u] k_u\rangle_{A^2(\Omega)}              \\
         & = \langle \psi, z k^\Deinf_u \rangle_{A^2(\Omega)}                           \\
         & = \langle \overline{z} \psi \mathds 1_\Omega, k^\Deinf_u \rangle_{A^2(\Deinf)} \\
         & = P_{A^2(\Deinf)}(\overline z \psi \mathds 1_{\Omega})(u).
    \end{align*}
    The function $z\mapsto \overline z \psi(z) \mathds 1_{\Omega}(z)$ is in $L^2(\Deinf)$ because $\psi \in A^2(\Omega)$ and $\Omega$ is bounded. We deduce that $f_0' \in A^2(\Deinf)$.
\end{delayedproof}

\begin{delayedproof}{th-analy-non-nc}
    Assume that $f_0\in \NCH(\omega)$. By \Cref{th-NCH-holom-necessary}, $f_0$ extends holomorphically to $\Cinf \setminus \cl{\omega}$ (\Cref{th-NCH-holom-necessary}). Moreover, denoting again this extension by $f_0$, we have $f_0(\infty) = 0$.

    Now, assuming that $f_0$ also extends holomorphically to $D(0,1+\epsilon)$, then $f_0$ is actually holomorphic on $\Cinf$. Since $f_0(\infty) = 0$, Liouville's theorem implies that $f_0 = 0$.
\end{delayedproof}

\begin{delayedproof}{th-CondNec}
    Let $f_0\in \NCH(\omega)$. According to \Cref{th-NCH-holom-necessary}, $f_0$ extends as a holomorphic functions of $\Cinf \setminus \cl \omega$ such that $f_0(\infty) = 0$ and $f_0 \in \Dirichlet(\Dext)$

    For a finite Borel measure $\mu$ on $\T$, we denote by $C\mu$ the \emph{Cauchy transform} of $\mu$ defined for $z\in \C\setminus \T$ by
    \[
        C\mu(z) = \int_\T \frac{\diff \mu(\zeta)}{1-\overline{\zeta}z}.
    \]

    \step{There exists $g\in L^2(\T)$ such that $f_0 = Cg$} (When writing $Cg$, we mean $C\mu$, where $\diff \mu(u) = g(u) \diff u/2\pi$.)
    Tumarkin's theorem~\cite[Theorem 5.3.1]{CMR06} \emph{almost} gives us this result. To actually get what we want, we find it is easier to adapt the proof of Tumarkin's theorem.

    Since $f_0$ extends analytically to $\Cinf\setminus \cl \omega$ we can write for $|u|<1$, $f_0(u) = \sum_{n\geq 0} a_n u^n$, and for $|u|>1$, ${f_0}(u) = \sum_{n<0} a_n u^n$ (there is no $a_0$ in this latter sum because ${f_0}(\infty) = 0$). Since ${f_0} \mathds 1_\D = f_0 \in H^2$, we have $a_n=\widehat{f_0}(n)$, $n\ge 0$, so that $(a_n)_{n\geq 0} \in \ell^2$. And since ${f_0} \mathds 1_\Deinf \in \Dirichlet$, we have $(\sqrt n
     a_n)_{n<0} \in \ell^2$. We can thus define $g$ a.e.\ on $\T$ as the $L^2(\T)$-function
    \[
        g(\zeta) = \sum_{n\geq 0}a_n \zeta^n - \sum_{n< 0} a_n \zeta^n,\quad \text{for almost every }|\zeta|=1.
    \]
    %Since $f_0$ extends analytically through $\T$ outside $\cl\omega$, we have $g(u)=0$ for almost every $$
    
    As in \eqref{CauchyTrans} we have
    \[
    \int_{\omega}\frac{g(w)}{1-\overline{w}z}\lvert\diff w|=
        \begin{cases}
            P_+g(z)=\sum_{n\ge 0}\widehat{g}(n)z^n=\sum_{n\ge 0}a_nz^n, & |z|<1,\\
            -P_-g(z)=-\sum_{n>0}\dfrac{\widehat{g}(-n)}{z^n}=\sum_{n>0}\dfrac{a_{-n}}{z^n}, & |z|>1,
        \end{cases}
    \]
    which shows that $f_0=Cg$ on $\Cinf\setminus \cl\omega$, and in particular $f_0=P_+g$ in $\D$. 

    \step{$\supp(g) \subset \cl\omega$}
    According to Fatou's jump theorem~\cite[Corollary 2.4.2]{CMR06}, for almost every $u\in \T$,
    \[
        g(u) = \lim_{r\to 1^-} (Cg(ru) - Cg(u/r)) = \lim_{r\to 1^-} ({f_0}(ru) - {f_0}(u/r)).
    \]
    But ${f_0}$ is holomorphic outside $\cl \omega$. Hence, if $u\in \T \setminus \cl \omega$, the above limit is $0$, hence $\supp g \subset \cl\omega$.
\end{delayedproof}

We next prove \cref{th-CondSuff}. This will be a consequence of a necessary and sufficient condition related to pseudo-Carleson measures (see \cref{CNS} below). We will first need some notation. 
We refer to \cref{tildeOm} for the following construction. 
    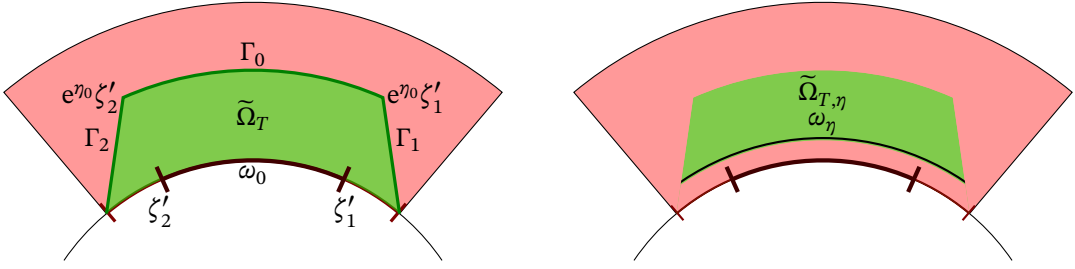
\begin{figure}
            \begin{center}
                \tikzmath{
    \q=1.1; \Rad=1.7; \rinz= 1.4; \rin = 1.1; \a=50; \b=130;  \shrink=0.2; \extend = 0.2;
    \az = (1-\shrink)*\a + \shrink*\b; \bz = \shrink*\a + (1-\shrink)*\b; 
    \ae = -\extend*\a + (1+\extend)*\b; \be = -\extend*\b + (1+\extend)*\a; 
    \abmid = (\a+\b)/2; \rinzmid = (1+\rinz)/2; \rinetamid = (\rin+2*\rinz)/3;}
% \tikzmath{
%     \q=1.1; \Rad=1.7; \rinz= 1.4; \rin = 1.1; \a=0; \b=180;  \shrink=0.4; \extend = 0.2;
%     \az = (1-\shrink)*\a + \shrink*\b; \bz = \shrink*\a + (1-\shrink)*\b; 
%     \ae = -\extend*\a + (1+\extend)*\b; \be = -\extend*\b + (1+\extend)*\a; 
%     \abmid = (\a+\b)/2; \rinzmid = (1+\rinz)/2; \rinetamid = (\rin+2*\rinz)/3;}
     
\begin{tikzpicture}[scale=3, every node/.style={fill opacity = 0.3, inner sep=1pt, text=black, text opacity=1}]
    \begin{scope}
        \draw (\ae:1) arc[start angle = \ae, end angle=\be, radius=1];
        
        \draw[fill opacity = 0.5, fill=Omegac] (\a:1) -- (\a:\Rad)
            arc[start angle=\a, end angle=\b, radius=\Rad]
            -- (\b:1)
            arc[start angle=\b, end angle=\a, radius=1];
        
        \draw[|-|, very thick, omegac] (\a:1) arc[start angle=\a, end angle = \b, radius = 1];%node[pos = 0.5, below left,fill=white, fill opacity = 0.3, text=black, text opacity = 1, inner sep = 1pt]{$\omega$};
        
        \fill[green, fill opacity=0.5](\a:1)
        -- (\az:\rinz)
        arc[start angle=\az, end angle=\bz, radius=\rinz]
        -- (\b:1)
        arc[start angle=\b, end angle=\a, radius=1];
    
        \draw[very thick, green!50!black] (\a:1)
        -- (\az:\rinz)
        node[pos=0.5, auto=right]{$\Gamma_1$}
        node[pos=1, right]{$\eu^{\eta_0} \zeta_1'$}
        arc[start angle=\az, end angle=\bz, radius=\rinz]
        node[pos=0.5, auto=right]{$\Gamma_0$}
        -- (\b:1)
        node[pos=0.5, auto=right] {$\Gamma_2$}
        node[pos=0, left]{$\eu^{\eta_0} \zeta_2'$};
    
        \draw[|-|, ultra thick, omegac!50!black] 
        (\az:1)  
        arc[start angle=\az, end angle = \bz, radius = 1] 
        node[pos = 0.5, auto=left]{$\omega_0$} 
        node[below, pos = 1, inner sep=5pt]{$\zeta'_2$}
        node[below, pos = 0, inner sep=5pt]{$\zeta'_1$};
    
        \node at (\abmid:\rinzmid){$\widetilde\Omega_{T}$};
    \end{scope}
    
    \begin{scope}[shift={(2.5,0)}]
        \draw (\ae:1) arc[start angle = \ae, end angle=\be, radius=1];
        
        \draw[fill opacity = 0.5, fill=Omegac] (\a:1) -- (\a:\Rad)
            arc[start angle=\a, end angle=\b, radius=\Rad]
            -- (\b:1)
            arc[start angle=\b, end angle=\a, radius=1];
        
        \draw[|-|, thick, omegac] (\a:1) arc[start angle=\a, end angle = \b, radius = 1];%node[pos = 0.5, below left,fill=white, fill opacity = 0.3, text=black, text opacity = 1, inner sep = 1pt]{$\omega$};

        \begin{scope}
            \path[clip] (\a:1)
            -- (\az:\rinz)
            arc[start angle=\az, end angle=\bz, radius=\rinz]
            -- (\b:1)
            arc[start angle=\b, end angle=\a, radius=1];

            \fill[green, fill opacity=0.5, even odd rule](\a:1)
            -- (\az:\rinz)
            arc[start angle=\az, end angle=\bz, radius=\rinz]
            -- (\b:1)
            arc[start angle=\b, end angle=\a, radius=1]
            (0,0) circle[radius=\rin-0.01];
            \draw[thick] (0,0) circle[radius=\rin];
        \end{scope}
        
        \node[above] at (\abmid:\rin){$\omega_\eta$};
    
        \draw[|-|, ultra thick, omegac!50!black] 
        (\az:1)  arc[start angle=\az, end angle = \bz, radius = 1];
    
        \node at (\abmid:\rinetamid){$\widetilde\Omega_{T,\eta}$};
    \end{scope}
    
\end{tikzpicture}
            \end{center}
            \caption{Left figure: illustration of $\widetilde\Omega_T$ (in green), with $\Omega_T$ in red, the paths $\Gamma_i$ ($i=0,1,2$) and the points $\zeta_i'$ and $\eu^{\eta_0}\zeta'_i$ ($i=1,2$). Right figure: illustration of $\widetilde \Omega_{T,\eta}$ (in green) and the path $\omega_\eta$.}
            \label{tildeOm}
    \end{figure}
Let $0<\eta<\eta_0<T$ and pick $\omega_0$ an arc compactly contained in $\omega$ with endpoints $\zeta_1'$ and $\zeta_2'$ having same center as $\omega$.
In order to avoid pathological situations, we will assume that $\zeta_k'$ is closely enough to $\zeta_k$ (for instance $\eta_0\pi/8\le |\zeta_k'-\zeta_k|\le \eta_0\pi/4$).  Introduce the segments $\Gamma_k=[\zeta_k,e^{\eta_0}\zeta'_k]$, $k=1,2$, and define  $\Gamma_{k,\eta}$ to be the part of $\Gamma_k$ in $|z|\ge (1+\eta)$. We also let $\omega_{\eta}$ to be the arc in $(1+\eta)\T$ comprised between the segments $\Gamma_1$ and $\Gamma_2$. Define $\widetilde{\Omega}_{T,\eta}$ as the set bounded by $\omega_{\eta}$, $\Gamma_{1,\eta}$, $\Gamma_0\coloneqq\omega_{\eta_0}$ and $\Gamma_{2,\eta}$. In particular
    \[
        \partial\widetilde{\Omega}_{T,\eta}=\omega_{\eta}\cup \Gamma_{1,\eta}\cup\Gamma_0\cup\Gamma_{2,\eta}.
    \]
    In view of a future application of Stokes' theorem we note that $\widetilde{\Omega}_{T,\eta}$ has regular boundary (in the sense of \cite[D\'efinition 2.3.1]{AM04}). Observe also that $\Gamma_0$ is at a uniform distance from $\partial\Omega^*_{T}$.
    %Finally, $\widetilde{\Omega}_{T}=\widetilde{\Omega}_{T,1}\cup \widetilde{\Omega}_{T,2}$. 

    We also let $\widetilde{\Omega}_{T}=\bigcup_{0<\eta<\eta_0}\widetilde{\Omega}_{T,\eta}$. Hence
    \[
        \partial\widetilde{\Omega}_{T}=\omega\cup \Gamma_1\cup\Gamma_0\cup\Gamma_2.
    \]

\begin{theorem}\label{CNS}
   The function $f_0\in H^2$ is in $\NCH(\omega)$ if and only if there exists $g\in L^2(\T)$, such that 
    \begin{enumerate}
        \item $f_0=P_+g$
        \item $\supp g\subset \cl\omega$
        \item $\sum\limits_{n<0} |n| |\widehat g(n)|^2 < +\infty$
        \item Denoting by $G$ the harmonic extension of $\zeta\to  \zeta g(\zeta)$ to $\D_e$ the measure $\overline{G}\diff \zeta$ on $\Gamma_{1}\cup\Gamma_{2}$ is pseudo-Carleson for $A^2(\Omega_T)$, i.e.
        \begin{equation}\label{PseudoCarl}
         \left|\int_{\Gamma_{1}\cup\Gamma_{2}}\overline{G}(\zeta)p(\zeta)\diff\zeta\right|
        \lesssim\|p\|_{A^2({\Omega}_{T})}
        \end{equation}
    \end{enumerate}

\end{theorem}

As we will see in the proof, $G(z)=g_1(1/\overline z)+g_2(z)$, $|z|>1$, with $g_2\in \Dirichlet(\Dext)$ and $g_1\in H^2$. In particular, $G$ can be written as $G(z) = \sum_{n\geq 0} G_n \overline z^{-n} + \sum_{n<0} G_n z^{-|n|}$ for some $(G_n)_{n\in\Z} \in \ell^2(\Z)$. Let us show that this implies that the integral on the left hand side of \eqref{PseudoCarl} is well defined. By rotation, we can assume without loss of generality that $\Gamma_i$ is a path of the form $s\in (0,\epsilon) \mapsto 1+s \eu^{\iu \theta}$, for some $\theta \in (-\pi/2,\pi/2)$. Thanks to Fubini's theorem,
\begin{align*}
    \Big|\int_{\Gamma_i} \overline{G(\zeta)} p(\zeta) \diff \zeta \Big|
    &= \Big| \int_{\Gamma_i} \left( \sum_{n\geq 0} \overline{G_n} \zeta^{-n} +\sum_{n<0} \overline{G_n} \overline \zeta^{-|n|}\right) p(\zeta) \diff \zeta\Big| \\
    & \lesssim \sum_{n\in \Z} \int_0^\epsilon |G_n| |1+s\eu^{\iu\theta}|^{-|n|} \|p\|_{\infty,\D}\diff s\\
    & \lesssim \sum_{n\in \Z} \int_0^\epsilon |G_n| (1+c_\theta s)^{-|n|} \|p\|_{\infty,\D}\diff s\\
    &\lesssim \|p\|_{\infty,\D} \sum_{n\in \Z} \dfrac{|G_n|}{|n|+1}\\
    &\lesssim \|p\|_{\infty,\D} \Big(\sum_{n\in \Z} \frac1{(|n|+1)^2}\Big)^{1/2} \Big(\sum_{n\in \Z} |G_n|^2\Big)^{1/2}.
\end{align*}

In particular, we have
\[
\int_{\Gamma_{1}\cup\Gamma_{2}}\overline{G}(\zeta)p(\zeta)\diff\zeta
=\lim_{\eta\to 0}\int_{\Gamma_{1,\eta}\cup\Gamma_{2,\eta}}\overline{G}(\zeta)p(\zeta)\diff\zeta.
\]

\begin{proof}
    We already know from \Cref{th-CondNec} that the conditions (i)-(iii) are necessary. In particular $g\in H^2+\mathcal{D}(\Dext)$ (defined initially on $\T$), i.e. %as in the theorem with 
    %\begin{equation}\label{BoundCond}
    %  \int_{\T}\frac{|g(e^{it})|^2}{|(e^{it}-\zeta_1)(e^{it}-\zeta_2)|}dt<\infty.
    %\end{equation}
%    By assumption, we have 
$g=g_1+g_2$ with $g_1\in H^2$ and $g_2\in \mathcal{D}(\Dext)$ (observe that this means that $-g_2$ is the analytic continuation of $g_1$ to $\Dext$ through $\T\setminus \cl\omega$). In our argument below we need to consider $zg$ rather than $g$ (this will be needed in the change of variable $\diff z=\iu z\diff t$ for $z=c\eu^{\iu t}$, where $c$ is a suitable constant), but clearly this does not affect our decomposition and we can actually write $\zeta g=g_1+g_2$ with $g_1\in H^2$ and $g_2\in \mathcal{D}(\Dext)$. 

\step{Use of Stokes' formula} Let $P^e[h]$ be the harmonic extension (Poisson integral) of a function $h\in L^2(\T)$ to the outer disk $\Dext$ (see for instance \cite[Theorem 4.11]{ABR01}):
\[
P^e[h](z)=\int_{\T}P(z,\eu^{\iu t})h(\zeta)\lvert\diff\zeta|,\quad |z|>1,
\]
where
\[
 P^e(z,\eu^{\iu t})=\frac{|z|^2-1}{|z-\zeta|^2},\quad |z|>1.
\]
Note that for a function $h$ in $H^2(\Dext)$ ($\simeq H^2_-=L^2(\T)\ominus H^2(\T)$) we have $h=P^e[h]$ where we identify $h$ defined in $\Dext$ with its non-tangential boundary values on $\T$. In particular, since $\mathcal{D}(\Dext)\subset H^2(\Dext)$, we have $P^e[g_2]=g_2$. Moreover, since $g_1\in H^2$ we have $h\coloneqq\overline{g}_1-\overline{g}_1(0)\in H^2_-\simeq H^2(\Dext)$ , which again can be extended to itself {\it via} the Poisson kernel: $P^e[h]=h$. Since the Poisson integral will also reproduce constants, and conjugating the expression, we obtain an anti-analytic function $P^e[g_1](z)=g_1(1/\overline{z})$. Hence
    \[
        P^e[zg](z)=g_2(z)+g_1\left(\frac{1}{\overline{z}}\right),\quad |z|>1.
    \]

    Let us now apply Stokes' formula to the form $w=\overline{G}p\diff z$ on the domain $\widetilde{\Omega}_{T,\eta}$,
    where $G=P^e[zg]$.
    Since $\widetilde{\Omega}_{T,\eta}$ is compactly contained in $\Dext$, all the functions $g, g_1, g_2$ are $C^{\infty}$ on this domain and we can apply Stokes' formula there.
    \begin{equation}\label{Stokes}
        \int_{\partial\widetilde{\Omega}_{T,\eta}}w(z)=
        \int_{\widetilde{\Omega}_{T,\eta}}\diff w(z)=
        \int_{\widetilde{\Omega}_{T,\eta}}\dzb(\overline{P^e[zg]}p)\diff \overline{z}\wedge \diff z
        =2\iu\int_{\widetilde{\Omega}_{T,\eta}}\overline{g_2'} p\diff A(z).
    \end{equation}
    Observe that the
    anti-analytic function $g_1(1/\overline{z})$ disappears when differentiating in $z$.
    
   Decomposing $\partial\widetilde{\Omega}_{T,\eta}$ this gives
    \[
       \int_{\partial\widetilde{\Omega}_{T,\eta}}w(z)= \int_{\omega_{\eta}} \overline{G}(z)p(z)\diff z
        +\int_{\Gamma_{1,\eta}\cup\Gamma_{2,\eta}}\overline{G}(z)p(z){\diff z}+\int_{\Gamma_0}\overline{G}(z)p(z){\diff z}=\int_{\widetilde{\Omega}_{T,\eta}}\overline{g_2'}(z) p(z)\diff A(z).
    \]

\step{Convergence of the integral on $\omega_{\eta}$, $\eta\to 0$} Let us show that
    \begin{equation}\label{Stokes-2}
        \int_{\partial\widetilde{\Omega}_{T}}\overline{G}(z) p(z) \diff z =
        2\iu\int_{\widetilde{\Omega}_{T}}\overline{g_2'}(z) p(z)\diff A(z).
    \end{equation}
    
    To see this, we start with Fubini's theorem:
    \[
    \int_{\omega_{\eta}} \overline{G(z)}p(z)\diff z
    =\int_{\omega_{\eta}} \int_{\T}\frac{|z|^2-1}{|z-\zeta|^2}\overline{\zeta g(\zeta)}\lvert\diff\zeta|
    p(z)\diff z
    =\int_{\T} \overline{\zeta g(\zeta)} \int_{\omega_{\eta}} \frac{|z|^2-1}{|z-\zeta|^2}p(z)\diff z
    \lvert\diff\zeta|.
    \]
    Let now $\tilde{\omega}_{\eta}$ be the projection of $\omega_{\eta}$ onto $\T$ (this is actually the intersection of the convex hull of $0$ and $\omega_\eta$ with $\T$). Then, with the change of variable $z=(1+\eta)\eu^{\iu t}$, so that $\diff z=(1+\eta)\iu\eu^{\iu t}\diff t=\iu z\diff t$, and $\zeta=\eu^{\iu s}$, we get
    \[
    \varphi_{\eta}(\zeta)\coloneqq\int_{\omega_{\eta}} \frac{|z|^2-1}{|z-\zeta|^2}p(z)\diff z
    =\int_{\tilde{\omega}_{\eta}}\frac{(1+\eta)^2-1}{|(1+\eta)\eu^{\iu (t-s)}-1|^2}(1+\eta)\iu\eu^{\iu t}    p((1+\eta)\eu^{\iu t})\diff t \to\begin{cases} \iu\zeta p(\zeta) & \zeta\in \omega\\
      0 &\zeta\notin \cl\omega\end{cases}
    \]
Indeed, it is clear that if $\zeta=\eu^{\iu s}\notin \cl\omega$, then the Poisson kernel goes uniformly to 0 on $\tilde{\omega}_{\eta}$. When $\zeta=\eu^{\iu s}\in\omega$, then for a sufficiently small $\eta>0$, we have $\zeta\in\tilde{\omega}_{\eta}$. Since the Poisson kernel is an approximate identity, it will concentrate in a neighborhood of $\zeta$ when $\eta\to 0$, so that in this case, the integral will converge to $\iu\zeta p(\zeta)$ (note that this function is a polynomial). Morever, 
\[
 |\varphi_{\eta}(\zeta)|=\Big\lvert\int_{\omega_{\eta}} \frac{|z|^2-1}{|z-\zeta|^2}p(z)\diff z\Big\rvert
 \lesssim \sup_{z\in \omega_{\eta}}|p(z)|\int_{\omega_{\eta}} 
 \frac{|z|^2-1}{|z-\zeta|^2}\lvert\diff z|
 \le \|p\|_{\infty,\T} \int_{\T}\frac{|z|^2-1}{|z-\zeta|^2}\lvert\diff z|
 =\|p\|_{\infty,\T},
\]
so that Lebesgue's dominated convergence theorem yields
    \begin{equation}\label{eq-Stokes-CVD}
        \lim_{\eta\to 0}
        \int_{\omega_{\eta}} \overline{G(z)}p(z)\diff z=\int_{\omega} \overline{g}(\zeta)p(\zeta)\lvert\diff \zeta|.
    \end{equation}

Together with Stokes' formula~\eqref{Stokes}, we get the desired formula~\cref{Stokes-2}.

    \step{Necessary condition}
    The preceding equality yields
    \begin{equation}\label{NC}
    \Big\lvert\int_{\Gamma_{1}\cup\Gamma_{2}}\overline{G}(z)p(z)\diff z\Big\rvert\le
        %\frac{1}{(1+\eta)^2}
        \Big|\int_{\omega} \overline{G}(z)p(z)\diff z\Big|+\Big|\int_{\Gamma_0} \overline{G}(z)p(z)\diff z\Big|+\Big|\int_{\widetilde{\Omega}_{T}}\overline{g_2'} (z)p(z)\diff A(z)\Big|
    \end{equation}

    Since $g_2\in \Dirichlet(\Dext)$, the Cauchy-Schwarz inequality yields
    \[
        \Big|\int_{\widetilde{\Omega}_{T}}\overline{g_2'}(z) p(z)\diff A(z)\Big|
        \le \|g_{2}\|_{\Dirichlet(\Dext)}\|p\|_{A^2(\Omega_{T})}.
    \]

    Also, since $\Gamma_0$ is compactly contained in $\Omega_{T}$ and hence $g$ is uniformly bounded on $\Gamma_0$, the integral on $\Gamma_0$ is controlled by $\|p\|_{A^2(\Omega_{T})}$ (using subharmonicity, it is easily seen that the norm of the point evaluation in $A^2(\Omega)$ is essentially controlled by the inverse of the distance to the boundary).

    By the assumption that $f_0$ is null-controllable, the modulus of the integral on $\omega$ in \cref{NC} is controlled by $\|p\|_{A^2(\Omega_T)}$. Hence, \cref{NC} shows that for every polynomial $p$,
    \[
        \Big|\int_{\Gamma_{1}\cup\Gamma_{2}}\overline{G}(z)p(z)\diff z\Big|
        \lesssim\|p\|_{A^2({\Omega}_{T})}
    \]
%    We deduce that $\overline{\zeta g}(\zeta)d\zeta|\Gamma_{1,\eta}\cup\Gamma_{2,\eta}$ is a pseudo-Carleson measure in the sense of \cite{Xiao00}. 
    %These are rather involved objects for which useful characterizations are not known to our knowledge

    \step{Sufficient condition}
    Suppose $f_0$ satisfies the 4 conditions of the theorem. Stokes' formula (\cref{Stokes-2}) yields
    \[
        \Big|\int_{\omega_{\eta}} \overline{G}(z)p(z)\diff z \Big|
        \le \Big|\int_{\Gamma_{1,\eta}\cup\Gamma_{2,\eta}}\overline{G}(z)p(z)\diff z\Big|
        %\frac{1}{(1+\eta)^2}
        +\Big|\int_{\Gamma_0} \overline{G}(z)p(z)\diff z\Big|+\Big|\int_{\widetilde{\Omega}_{T,\eta}}\overline{g_2'} p\diff A(z)\Big|.
    \]

    The last two terms are again controlled by $\|p\|_{A^2(\Omega_T)}$ using the same arguments as in the necessary condition (note that $g\in H^2+\mathcal{D}(\Dext)$ is our standing hypothesis). Moreover, we have by assumption $|\int_{\Gamma_{1}\cup\Gamma_{2}}\overline{G}(z)p(z)\diff z|
\lesssim \|p\|_{A^2(\Omega_T)}$, from where we deduce that
\[
\Big|\int_{\omega} \overline{g}(\eu^{\iu t})p(\eu^{\iu t})\diff t\Big|\lesssim  \|p\|_{A^2(\Omega_T)}
\]
%Clearly, the operator $R_{\eta}$ initially defined for polynomials extends boundedly to $A^2(\Omega_T)$ for every $0<\eta<\eta_0$. Moreover, $\lim_{\eta\to 0} R_{\eta}p=:R_0 p=\int_{\Gamma_1\cup \Gamma_2}\overline{G}(\zeta)p(\zeta)\diff\zeta$ is well defined, and by assumption $R_0$ extends boundedly to $A^2(\Omega_T)$. We deduce that the orbits $R_{\eta}f$ are bounded for $f\in A^2(\Omega_T)$, so that again by the uniform boundedness principle we can show that the norm of $T_{\eta}$ is uniformly bounded. As in step 1 -- which only uses the fact that $G$ is the harmonic extension of an $L^2$-function and $p$ is a polynomial -- we see that  
%\[
%|\int_{\omega} \overline{g}(\eu^{\iu t})p(\eu^{\iu t})\diff t|=
%\lim_{\eta\to 0}|\int_{\omega_{\eta}} \overline{G(z)}p(z)\diff z|\le C\|p\|_{A^2(\Omega_T)}, 
%\]
showing that $g\in \NCH(\omega)$.
\end{proof}

We are now in a position to prove Theorem \ref{th-CondSuff}.

\begin{delayedproof}{th-CondSuff}
%   It is enough to how the condition
%    \[
%     \int_{\T}\frac{|f(\eu^{\iu t})|^2}{|(\eu^{\iu t}-\zeta_1)(\eu^{\iu t}-\zeta_2)|}\diff t<\infty
%    \]
%    implies the pseudo-Carleson measure condition. 
For this, it is enough to estimate the integral condition of \Cref{CNS} on $\Gamma_{1,\eta}\cup\Gamma_{2,\eta}$ under the assumption
    \begin{equation}\label{BoundCond}
        \int_{\T}\frac{|g(\eu^{\iu t})|^2}{|(\eu^{\iu t}-\zeta_1)(\eu^{\iu t}-\zeta_2)|}\diff t<\infty.
    \end{equation}
    %Let us take care of the the integral on $\Gamma_{k,\eta}$, 

     It is clearly enough to check that
    \[
        \int_{\Gamma_1}|\overline{G}(z)p(z)|\lvert\diff z|\lesssim \|p\|_{A^2(\Omega^*_{T})}.
    \]
    %By triangular inequality, this is certainly true if the measure $d\mu(z)=g(z) d_{\Gamma_1}(z)$, where $d_{\Gamma_1}(z)$ is arc-length measure on $\Gamma_1$, is a Carleson measure for $A^1(\Omega)$, since then for every polynomial
    %\[
    % |\int_{\Gamma_1}|\overline{g(z)}p(z)||dz||\lesssim \|p\|_{A^1(\Omega^*_{T})}\lesssim \|p\|_{A^2(\Omega^*_{T})}.
    %\]
    %With the Poisson kernel for $\Dext$ in mind, we get
    %\[
    % g(z)=\frac{1}{2\pi}\int_{\T}\frac{|z|^2-1}{|z-e^{it}|^2}g(e^{it})dt,\quad |z|>1.
    %\]
    Let us parameterize $\Gamma_1$. In order to simplify the argument we assume that $\omega$ starts at $\zeta=1$ and that the segment $\Gamma_1$ takes the direction $(1+i)$. Then a parametrization is given by $\gamma(s)=1+s(1+i)$, $s\in [0,\delta_0]$ (we assume that $\zeta_1$ is the lower left corner of $\Omega_T$, and so $\Gamma_1$ moves north-east). In particular $\lvert\diff z|=\sqrt{2}\diff s$. By the Cauchy-Schwarz inequality, we have
    \[
        \int_{\Gamma_1}|\overline{G}(z)p(z)|\lvert\diff z|\lesssim \left(\int_0^{\delta_0}\frac{|G(\gamma(s))|^2}{s}\diff s\right)^{1/2} \left(\int_0^{\delta_0} s|p(\gamma(s)|^2 \diff s\right)^{1/2}
    \]
    %The first integral on the right hand side is controlled by our %assumption \eqref{BoundCond}. 
    We claim that the second term translates to a Carleson measure result. Indeed, setting $\diff\mu(z)=s \diff_{\Gamma_1}(\gamma(s))$, where $\diff_{\Gamma_1}$ is arclength measure on $\Gamma_1$, it is enough to check that $\diff\mu$ is a Carleson measure for $A^2(\Omega_{T})$. For this, using a result by Gonzalez \cite{Gonzalez20}, we have to check that $\mu(B(z,r))\lesssim r^2$ for Whitney balls in $\Omega_{T}$, which are balls having a radius which is less than one half the distance %$\delta_{\Omega^*_{T}}$ 
    of $z$ to the boundary of $\Omega_{T}$ but still is comparable to that distance. Clearly, in the worst case such a ball contains $\Gamma_1$ as a diameter, in which case $\int_{B(z,r)\cap\Gamma_1}s \diff s\simeq r^2$.

    In order to finish the proof we thus have to estimate the first term. By Jensen's inequality, we have
    \[
        |G(\gamma(s))|^2=\left(\int_{\T}\frac{|\gamma(s)|^2-1}{|\gamma(s)-\eu^{\iu t}|^2}g(\eu^{\iu t})\diff t\right)^2\lesssim \int_{\T}\frac{|\gamma(s)|^2-1}{|\gamma(s)-e^{\iu t}|^2}|g\eu^{\iu t})|^2\diff t.
    \]
    Now, letting $L$ be the length of $\omega$ (which we had chosen to start at $\zeta_1=1$), and  
since for $z=\gamma(s)\in\Gamma_1$ we have $|z|^2-1\simeq s$ and $|\eu^{\iu t}-z|^2\simeq  (s-t)^2+s^2$, an application of Fubini's theorem gives 
    \begin{equation}\label{PoissonEstimate}
        \int_{0}^{\delta_0} \frac{|G(\gamma(s))|^2}{s}\diff s
        \lesssim \int_0^L |g(\eu^{\iu t})|^2\int_0^{\delta_0}\frac{1}{(s-t)^2+s^2}\diff s\diff t.
    \end{equation}
Observing that $\int_0^{\delta_0}\frac{1}{(s-t)^2+s^2}\diff s=\frac{1}{2}\int_0^{\delta_0}\frac{1}{(s-t/2)^2+3t^2/4}\diff s$, by a change of variables and standard estimates we control this integral by $C/t$ (starting from the last integral, one can reach this estimate also by reinterpreting this integral, after multiplication by $t/2$, as a Poisson integral).
By \eqref{BoundCond} we obtain the required control on the piece $\Gamma_1$, and similarly on $\Gamma_2$.
\end{delayedproof}

To prove that $\NCH(\omega)$ is dense, we will use the fact that $P_+$ is self-adjoint on $W^{k,2}(\T)$. Recall also that $H^2\cap W^{k,2}(\T)$ is endowed with the $W^{k,2}$ scalar product given by
\begin{equation*}
    \langle f,g\rangle_{W^{k,2}(\T)} \coloneqq \sum_{n\in\Z} (1+n^2)^k \widehat f(n) \overline{\widehat g(n)}.
\end{equation*}
We will also use the following
\begin{lemma}\label{th-lemma-ortho}
Let $k\in \N$. Let $f\in H^2\cap W^{k,2}(\T)$. Assume that
\begin{equation*}
    \forall g\in C_c^\infty(\omega),\ \langle f,g\rangle_{W^{k,2}(\T)}=0.
\end{equation*}
Then $f=0$.
\end{lemma}
\begin{proof}
    In this proof, it is useful to distinguish the maps $z\in \D \mapsto f(z)$ and $\widetilde f \colon x\in \R \mapsto f(\eu^{\iu x})$. In the same vein, we define $\widetilde \omega\subset [-\pi,\pi]$ by $\omega = \{\eu^{\iu x},\ x\in\widetilde\omega\}$. Up to a rotation, we may assume that $-1\notin \cl \omega$, in which case $\widetilde \omega$ is an interval of $[-\pi,\pi]$.
    
    \step{$\widetilde f(x)$ is a linear combination of $x^\ell \eu^{\pm x}$ on $\omega$}
    Indeed, if $g\in C_c^\infty(\omega)$, the $W^{k,2}$ scalar product can be written as
    \begin{equation*}
        \langle f,g\rangle_{W^{k,2}(\T)} = \langle (I-\partial_x^2)^k \widetilde f,\widetilde g\rangle_{W^{-k,2}(-\pi,\pi),W^{k,2}_0(-\pi,\pi)}.
    \end{equation*}
    If this is $0$ for every $g\in C^\infty_c(\omega)$ (i.e., for every $\widetilde g\in C_c^\infty(\widetilde \omega)$), this means that $(I-\partial_x^2)^k \widetilde f = 0$ on $\widetilde \omega$ (by definition in the distribution sense). The functions $x\mapsto x^\ell \eu^{\pm x}$ for $0\leq\ell < k$ are a basis of solutions of this ODE, hence there exists $(a^\pm_\ell)_\ell$ such that $\forall x\in\widetilde \omega$,
    \begin{equation}\label{eq-f-omega}
        \widetilde f(x) = 
        \sum_{0\leq \ell<k}  \iu^\ell x^\ell(a^+_\ell\eu^{-x} + a^-_\ell \eu^{x}).
    \end{equation}
    The $\iu^\ell$ is only here to simplify the notations in the following.

    \step{Extension of the formula~\eqref{eq-f-omega} to a larger domain} 
    We can define the holomorphic function $h$ on $\C\setminus (-\infty,0]$ by
    \begin{equation}\label{eq-def-g-dense}
        h(z) \coloneqq 
        \sum_{0\leq \ell<k}  (\ln(z))^\ell(a^+_\ell z^\iu + a^-_\ell z^{-\iu}).
    \end{equation} 
    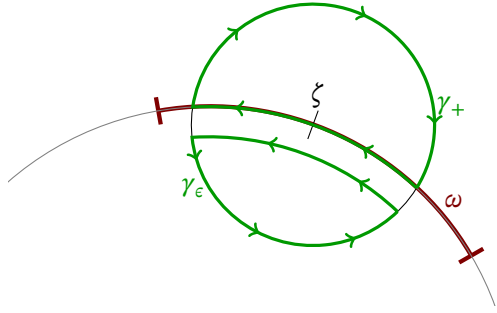
\begin{figure}
            \begin{center}
                \tikzmath{
    \rotation = 70; \omegastart = -40; \omegaend=30; \lx = 1; \ly = 0.6; \thr= 20; \eps = 0.1; \radius = 0.4; \tick = 0.05;
    \rin = 1- \eps;
    \xa = 1-\radius*\radius/2; \ya = sqrt(1-\xa*\xa); \thetaa = atan(\ya/\xa); \thetap = 180+atan(\ya/(\xa-1));
    \xb = (1+\rin*\rin-\radius*\radius)/2; \yb = sqrt(\rin*\rin-\xb*\xb); \thetab = atan(\yb/\xb); \thetabp = 180+atan(\yb/(\xb-1));}
     
\begin{tikzpicture}[scale=4, every node/.style={fill opacity = 0.3, inner sep=1pt, text=black, text opacity=1}]
    \path[clip] (\rotation:1) +(-\lx,-\ly) rectangle +(\lx,\ly);
    \begin{scope}[rotate = \rotation]
        \draw[|-|, ultra thick, omegac] (\omegastart:1) arc[start angle = \omegastart, end angle = \omegaend, radius = 1] node[pos=0.15, above right, text=omegac]{$\omega$};
        
        \draw (0:1)  +(\tick,0)node[above]{$\zeta$} -- +(-\tick,0) +(0,-\tick) -- +(0,\tick);
        
        \draw (0:1) circle[radius = \radius];
        
        \draw[thin, gray] (0,0) circle[radius = 1];

        \begin{scope}[decoration={markings,
            mark= between positions 0.1 and 0.9 step 0.2 with {\arrow{>}}}]
            \draw[very thick, green!60!black, postaction={decorate}] (-\thetaa:1) arc[start angle = -\thetaa, end angle = \thetaa, radius = 1] arc[start angle = \thetap, end angle = -\thetap, radius = \radius] node[pos=0.8,right, text=green!60!black]{$\gamma_+$}; 
    
            \draw[very thick, green!60!black, postaction={decorate}] (-\thetab:\rin) arc[start angle = -\thetab, end angle = \thetab, radius = \rin] arc[start angle = \thetabp, end angle = 360-\thetabp, radius = \radius] node[pos=0.2,left, text=green!60!black]{$\gamma_\epsilon$};
        \end{scope}
    \end{scope}

\end{tikzpicture}
            \end{center}
            \caption{Illustration of the paths $\gamma_\epsilon$ and $\gamma_+$. They are used in the proof of \cref{th-lemma-ortho} to prove with Morera's theorem that the function $H(z) = f(z)$ if $|z|<1$ and $f(z) = h(z)$ if $|z|\geq 1$ and $z\in\notin(-\infty,0]$ is holomorphic.}
            \label{fig-morera}
    \end{figure}
    By direct inspection we see that the non-tangential limits of $f$ on $\omega$ coincide with the values of $h$ on $\omega$. This implies that $f$ extends $h$ through $\omega$ inside $\D$. To see this, one can for instance use Morera's theorem. Define $H(z)=f(z)$ if $|z|<1$ and $H(z)=h(z)$ for $|z|\ge 1$, $z\notin (-\infty,0]$. For $\zeta\in\omega$, let $D_{\zeta}$ be a circle centered at $\zeta$ such that $D_{\zeta}\cap \T\subset \omega$, and define $\gamma=\partial D_{\zeta}$. For every $\epsilon>0$, we let $\gamma_{\epsilon}=(\gamma\cap (1-\epsilon)\cl\D)\cup ((1-\epsilon)\T\cap D_{\zeta})$ (see \cref{fig-morera}). Define also $\gamma_+=\partial (D_{\zeta}\cap \cl\Dext)$. Since $f\in \Hol(\D)$ and $h\in\Hol(\Dext)$, $\int_{\gamma_{\epsilon}}H(z)\diff z=\int_{\gamma_{\epsilon}}f(z)\diff z=0$, and $\int_{\gamma_+}H(z)\diff z=\int_{\gamma_+}g(z)\diff z=0$. Moreover, $f\in H^2$ implies that 
    \[
     \int_{(1-\epsilon)\T\cap D_{\zeta}}H(z)\diff z=\int_{(1-\epsilon)\T\cap D_{\zeta}}f(z)\diff z\xrightarrow[\epsilon\to 0]{} \int_{D_{\zeta}\cap\T}f(z)\diff z=\int_{D_{\zeta}\cap\T}H(z)\diff z.
    \]
    Hence 
    \[
     \int_\gamma H(z)\diff z=\lim_{\epsilon\to 0}
     \left(\int_{\gamma_{\epsilon}} H(z)\diff z+\int_{\gamma_+} H(z)\diff z \right)=0.
    \]
    As a result, $H$ is holmorphic in $\C\setminus (-\infty,0]$, and hence $f(z)=h(z)$, $z\in \D\setminus (-1,0]$.
    
    \step{Conclusion using the holomorphy of $f$}
    We have proved that $f=h$. But $f\in H^2$, and as such, is holomorphic on $\D$, while the definition of $h$~\eqref{eq-def-g-dense} involves functions with a cut along $\R_-$. The aim is to prove that all the coefficients in \eqref{eq-def-g-dense} are zero by using that for all $0<r<1$,  $h(-r+\iu 0) = h(-r-\iu 0)$.
    
    Assume that at least one of the $a_\ell^\pm\neq 0$ and set $\ell_0$ the largest $\ell$ such that $(a_\ell^+, a_\ell^-) \neq (0,0)$. For $r\to 0^+$,
    \[
        h(-r\pm\iu0) = \lim_{\epsilon\to 0^+}h(-r\pm \epsilon) = 
        \sum_{0\leq \ell < k} (\ln(r)\pm \iu\pi)^\ell \left(a_\ell^+ \eu^{\iu \ln(r) \mp \pi} + a_\ell^- \eu^{-\iu \ln(r) \pm\pi} \right).
    \]
    Setting $r = \eu^s$ and taking $s\to-\infty$, we write this as
    \begin{align*}
        h(-\eu^s\pm\iu0) 
        &= (s\pm\iu \pi)^{\ell_0}\left(a_{\ell_0}^+ \eu^{\iu s \mp \pi} + a_{\ell_0}^- \eu^{-\iu s +\pm\pi}\right) + \bigO(s^{{\ell_0}-1})\\
        &= s^{\ell_0}\left(a_{\ell_0}^+ \eu^{\iu s \mp \pi} + a_{\ell_0}^- \eu^{-\iu s +\pm\pi}\right) + \bigO(s^{{\ell_0}-1})
    \end{align*}
        
    But $f=h$ is supposed to be holomorphic on $\D$. Hence, for $s<0$, $h(-\eu^s +\iu 0) = h(-\eu^s - \iu 0)$. In particular,
    \begin{align*}
        0 &= h(-\eu^s + \iu 0) - h(-\eu^s - \iu 0)\\
        &= s^{\ell_0}2\sinh(\pi)\left(-a^+_{\ell_0} \eu^{\iu s} + a_{\ell_0}^- \eu^{-\iu s}\right) + \bigO(s^{\ell_0-1}).
    \end{align*}
    This is only possible if 
    \[
        -a^+_{\ell_0} \eu^{\iu s} + a_{\ell_0}^- \eu^{-\iu s} \xrightarrow[s\to -\infty]{} 0.
    \]
    This is in turn only possible if $a^+_{\ell_0} = a^-_{\ell_0}=0$, a contradiction with $(a^+_{\ell_0},a^-_{\ell_0}) \neq (0,0)$.
    
    By contradiction, all the $a^\pm_\ell$ are $0$, hence $f = h = 0$.
\end{proof}
With this proposition, we can easily prove that $\NCH(\omega)$ is dense:
\begin{delayedproof}{th-dense}
    Let $h\in (\NCH(\omega)\cap W^{k,2}(\T))^\bot$ (orthogonal in $H^2\cap W^{k,2}(\T)$). Every $g\in C_c^\infty(\omega)$ satisfies the hypothesis of our sufficient condition (\cref{th-CondSuff}) so that $P_+g\in \NCH(\omega)$, and hence % Hence, for every $g\in C_c^\infty(\omega)$,
    \begin{equation*}
        0 = \langle h,P_+g\rangle_{W^{k,2}(\T)}
          = \langle h,g\rangle_{W^{k,2}(\T)}.
    \end{equation*}
    Thus, according to \cref{th-lemma-ortho} $h = 0$. This proves that $\NCH(\omega)^\bot = \{0\}$.
\end{delayedproof}

Finally, let us prove the necessary and sufficient conditions for $W^{1/2,2}$ initial states.
\begin{delayedproof}{th-CNS-W12}
    \step{Direct implication}
    Let $f_0\in \NCH(\omega)\cap W^{1/2,2}(\T)$. According to the necessary condition of \cref{th-CondNec}, there exists $g\in L^2(\T)$ such that $f_0 = P_+g$, $\supp g\subset \cl\omega$ and $\sum_{n<0} |n| |\widehat g(n)|^2 < +\infty$. Since $f_0 = P_+g\in W^{1/2,2}(\T)$ by hypothesis, we also have $\sum_{n\geq 0} |n||\widehat g(n)|^2 < +\infty$, thus $g\in W^{1/2,2}(\T)$.
    
    \step{Reciprocal implication}
    Let $g\in W^{1/2,2}(\T)$ with $\supp g \subset \cl \omega$. Of course, $P_+g\in H^2 \cap W^{1/2,2}(\T)$. Let us prove $P_+g\in \NCH(\omega)$.
    
    By definition, $g\in W^{1/2,2}_{00}(\omega)$~\cite[Definition 1.85]{Leoni23}. Hence, $\int_{\T}\frac{|f(\eu^{\iu t})|^2}{|(\eu^{\iu t}-\zeta_1)(\eu^{\iu t}-\zeta_2)|}\diff t<\infty$ (see, e.g., \cite[Theorem 1.87]{Leoni23} for the analog property for $W^{1/2,2}_{00}(\R_+)$, and use some partition of unity).
    
    Thus, according to our sufficient condition (\cref{th-CondSuff}), $P_+ g \in \NCH(\omega)$.
\end{delayedproof}
The reciprocal implication can also be shown without \cref{th-CondSuff} in the following way. Consider $\widetilde g\colon \partial\Omega\to\C$ defined by $\widetilde g = g$ on $\omega$, $\widetilde g = 0$ on $\partial\Omega\setminus \omega$. Since $g\in W^{1/2,2}_{00}(\omega)$, $\widetilde g \in W^{1/2,2}(\partial\Omega)$. Hence, there exists $h\in W^{1,2}(\Omega)$ such that the trace of $h$ on $\partial\Omega$ is $\widetilde g$~\cite[18.40]{Leoni17}. Then, we can solve the Dirichlet problem (in the sense of distributions) $-\Delta u = \Delta h$ on $\Omega$, $u\in W^{1,2}_0(\Omega)$. Notice that $v = u+h$ is the solution of $\Delta v = 0$ on $\Omega$, $v = \widetilde g$ on $\partial\Omega$. Since $v\in W^{1,2}(\Omega)$ a use of Stokes formula will give the result (see \cref{sec-formal-stokes} for more details).

\subsection{Characterization thanks to the Poisson kernel}\label{sec-h2-poisson}
We will now prove the necessary and sufficient condition of \cref{th-CNS}. We start presenting the main idea before giving a rigorous proof.

\subsubsection{Formal argument}\label{sec-formal-stokes}
The main idea of the proof is essentially the same as in the proof of \cref{CNS} and based on Stokes' formula \eqref{Stokes}. Now, instead of extending $zg$ harmonically to $\Dext$ we use harmonic extension to $\Omega$. Let $g \in L^2(\T)$ with $\supp g\subset \cl\omega$. As in the rest of the article, we consider $\Omega = \Omega_T$ as defined in \cref{fig-Omega}. We denote by $u$ the solution of (see \cref{sec-conformal-and-poisson,th-poisson-lipschitz} for the proper definition):
\[
 \left\{
 \begin{IEEEeqnarraybox}[][c]{l"l}
      \Delta u(z) = 0,&  z\in\Omega\\
      u(z) = zg(z),& x \in \omega\\
      u(z) = 0,& z\in \partial\Omega\setminus\omega.
 \end{IEEEeqnarraybox}
 \right.
\]
Here and in the following, $zg$ is just the function $z\mapsto zg(z)$. Let $p\in \C[X]$. Since $P_+$ is self adjoint on $L^2(\T)$ and $p\in H^2$,
\begin{IEEEeqnarray*}{+rCl+s+}
    \langle p,P_+g\rangle_{H^2} &=&\langle p,g\rangle_{H^2}\\
    &=&\int_{\omega} p(z) \overline{g(z)} \lvert\diff z| & ($\supp g\subset \cl\omega$)\\
    &=&\int_{\omega} zp(z) \overline{zg(z)} \lvert\diff z| & ($|z|^2=1$ on $\omega$)\\
    &=& -\iu\int_{\omega} p(z) \overline{zg(z)}\diff z &(on $\omega$, $\diff z = \iu z \lvert\diff z|$)\\
    &=& -\iu\int_{\partial\Omega} p(z) \overline{u(z)}\diff z & (boundary values of $u$).
    \intertext{Consider the differential form $\omega=p(z) \overline u(z)\diff z$ (where $u$ is now the harmonic extension to $\Omega$ instead of $G=P^e[zg]$ appearing in \eqref{Stokes}), an application of Stokes formula gives formally:}
    \langle p,P_+g\rangle_{H^2}
    &=& -\iu\int_{\Omega} \dzb (p(z)\overline{u(z)}) \diff \overline z \wedge \diff z\\
    &=& 2\int_{\Omega} p(z)\overline{\dz u}(z) \diff A(z) &($\dzb (p)=0$).
\end{IEEEeqnarray*}
Hence, if we want to estimate the left-hand side by $\|p\|_{L^2(\Omega)}$, it should be necessary and sufficient to have $\dz u\in L^2(\Omega)$.

\subsubsection{Rigorous proof}
If we assume $\dz u\in L^2(\Omega)$, we can justify the computation above by approximating $\Omega$ by a smooth subdomain (see \cref{sec-conformal-and-poisson,th-weak-stokes}). Hence, the sufficient part of \cref{th-CNS} is proved.

However, the proof of the necessary part of \cref{th-CNS} will require additional ideas. The observability inequality (\cref{th-MainEstim}) allows us to define $\psi_g$ by:
\begin{definition}\label{def-psi-g}
    Let $g\in L^2(\T)$ with $\supp g\subset \cl\omega$. If $P_+g \in \NCH(\omega)$, we denote by $\psi_g$ the unique function in $A^2(\Omega)$ such that
    \[
    \forall p\in \C[X],\ \langle p,P_+g\rangle_{H^2} = \int_{\Omega} p(z) \overline{\psi_g(z)} \diff A(z).
    \]
\end{definition}

If $g$ is regular enough, we can explicitely identify $\psi_g$:
\begin{prop}\label{th-formula-psi_g}
    If $g\in C_c^\infty(\omega)$ with $P_+g \in \NCH(\omega)$, then $\psi_g = 2\dz P^{\Omega}(zg)$.
\end{prop}
\begin{proof}
    If $g\in C_c^\infty(\omega)$ and $g$ is extended by 0 on $\partial\Omega\setminus\omega$, standard elliptic regularity methods tell us that $P^{\Omega} g\in W^{1,2}(\Omega)$ (see \cref{th-elliptic-reg} for a proof relying on the Riemann map), and this is also true for $P^{\Omega} (zg)$ since $zg$ has the same regularity as $g$. Hence, $\dz P^{\Omega}(zg)\in L^2(\Omega)$ and we can apply the Stokes formula of \cref{th-weak-stokes}. As in \cref{sec-formal-stokes}, we get
    \begin{align*}
    \langle p, P_+ g\rangle_{H^2} &= \langle p, g\rangle_{H^2}\\
    &= 2\int_{\Omega} p(z) \overline{\dz P^{\Omega}(zg)(z)} \diff A(z).
    \end{align*}
    
    Finally, since $P^{\Omega}(zg)$ is harmonic and since $\Omega$ is simply connected, $P^{\Omega}(zg)$ can be written as a sum
    of a holomorphic and an anti-holomorphic function, say $P^{\Omega}(zg)=f_1 + f_2$. Since $f_2$ is anti-holomorphic,
    $\dz \overline f_2=0$. Hence, $\dz P^{\Omega}(zg)=f_1'$ is holomorphic. Since we have already proved that $\dz P^{\Omega}(zg)\in L^2(\Omega)$, it follows that $P^\Omega(zg)\in A^2(\Omega)$.
\end{proof}

Now, the necessary part of \cref{th-CNS} will be an easy consequence of the following two propositions:
\begin{prop}\label{th-cv-dzPg}
    Let $g,g_n \in L^2(\T)$ be supported in $\cl\omega$. Assume that $g_n \xrightarrow[n\to \infty]{} g$ in $L^2(\omega)$. Then, for every $z\in \Omega$,
    \[
    \dz P^{\Omega}(zg_n)(z) \xrightarrow[n\to+\infty]{} \dz P^{\Omega}(zg)(z).
    \]
\end{prop}

\begin{prop}\label{th-cv-psi}
    Let $g,g_n \in L^2(\T)$ be supported in $\cl\omega$ and such that $P_+g_n\in \NCH(\omega)$ and $P_+ g\in \NCH(\omega)$. Assume that $g_n \xrightarrow[n\to \infty]{} g$ in $L^2(\omega)$. Then, for every $z\in \Omega$,
    \[
    \psi_{g_n}(z) \xrightarrow[n\to+\infty]{} \psi_g(z).
    \]
\end{prop}

\begin{delayedproof}{th-cv-dzPg}
    The map $g\mapsto P^{\Omega}g$ is continuous from $L^2(\omega)$ to $h(\Omega)$ endowed with the topology of uniform convergence on compact subsets (\cref{th-poisson-lipschitz}). Since $P^{\Omega}g$ is harmonic, for every $z\in\Omega$ and every $r>0$ such that $\overline{D(z,r)}\subset \Omega$, a standard estimate on harmonic functions~\cite[Chapter 2, Theorem 7]{Evans10} tells us:
    \[
    |\dz P^{\Omega}g(z)|\leq C_r\|P^{\Omega}g\|_{L^1(D(z,r))} \leq C_r'\|P^{\Omega}g\|_{L^\infty(D(z,r))}.
    \]
    This proves that for every $z\in\Omega$, the map $g\in L^2(\omega)\mapsto \dz P^{\Omega}g(z)$ is continuous.
\end{delayedproof}

\begin{delayedproof}{th-cv-psi}
    \step{Representation of $\psi_g$} 
    We claim that for all $h\in L^2(\omega)$ such that $P_+h \in \NCH(\omega)$,
    \begin{equation}\label{eq-psi-g}
        \psi_h(z) = \int_\omega h(\zeta) \overline{k_z^{\Omega}(\zeta)}\lvert\diff\zeta|,
    \end{equation}
    where $k^{\Omega}_z$ is the Bergman kernel on $\Omega$. Recall that $k^\Omega_z$ is the unique function in $A^2(\Omega)$ such that for all $f\in A^2(\Omega)$,
    \[
        f(z) = \langle f,k_z^{\Omega}\rangle_{A^2(\Omega)}.
    \]
    
    Observe that using the boundary behaviour of the conformal map $\D\to\Omega$, \cref{th-berg-continuous} allows to show that for all $z\in\Omega$, $k^{\Omega}_z$ extends continuously on $\cl{\Omega}$. Hence, according to Mergelyan's theorem~\cite[Theorem 20.5]{Rudin86}, there exists a sequence of polynomials $(p_k)_k$ such that $p_k \to k^{\Omega}_z$ in $L^\infty(\Omega)$.
    
    Thus
    \begin{IEEEeqnarray*}{+rCl+s}
        \psi_h(z) 
        &=& \int_{\Omega}\psi_h(\zeta) \overline{k_z^{\Omega}(\zeta)} \diff A(\zeta) 
        & (definition of the Bergman kernel)\\
        &=& \lim_{k\to\infty} \int_{\Omega}\psi_h(\zeta) \overline{p_k(\zeta)} \diff A(\zeta) & ($p_k\xrightarrow[k\to\infty]{} k_z^{\Omega}$)\\
        &=& \lim_{k\to\infty} \langle P_+h,p_k\rangle_{H^2}& (definition \ref{def-psi-g} of $\psi_h$)\\
        &=& \lim_{k\to \infty} \langle h,p_k\rangle_{L^2} & ($P_+$ selfadjoint and $P_+p_k = p_k$)\\
        &=& \lim_{k\to \infty} \int_\omega h(\zeta)\overline{p_k(\zeta)}\lvert\diff\zeta| & ($\supp h\subset \cl\omega$)\\
        &=& \int_\omega h(\zeta) \overline{k_z^{\Omega}(\zeta)} \lvert\diff \zeta| &  ($p_k\xrightarrow[k\to\infty]{} k_z^{\Omega}$).
    \end{IEEEeqnarray*}
    This is the claimed formula.
    
    \step{Conclusion}
    We can apply this formula for $g$ and $g_n$. Thus,
    \begin{align*}
        \psi_g(z) 
        &= \int_\omega g(\zeta)\overline{k_z^\Omega(\zeta)}\lvert\diff\zeta|\\
        &= \lim_{n\to +\infty} \int_\omega g_n(\zeta)\overline{k_z^{\Omega}(\zeta)}\lvert\diff\zeta|\\
        &= \lim_{n\to\infty} \psi_{g_n}(z).\qedhere
    \end{align*}
\end{delayedproof}

Combining these results, we can now prove the following formula for $\psi_g$:
\begin{prop}
    Let $g\in L^2(\T)$, be such that $\supp g \subset \cl\omega$ and $P_+g\in\NCH(\omega)$. Then,
    \[
        \psi_g = 2\dz P^{\Omega}(zg).
    \]
\end{prop}
Since $\psi_g\in A^2(\Omega)$, this formula completes the proof of the necessary part of \cref{th-CNS}. 
\begin{proof}
    Let $g_n\in C_c^\infty(\omega)$ be such that $g_n\to g$ in $L^2(\T)$. According to \cref{th-CNS-W12}, $P_+g_n\in \NCH(\omega)$, and according to \cref{th-formula-psi_g}, $\psi_{g_n} = 2\dz P^{\Omega}(zg_n)$.
    
    According to \cref{th-cv-dzPg},
    \[
        \psi_{g_n} = 2\dz P^{\Omega}(zg_n) \xrightarrow[n\to\infty]{} 2\dz P^{\Omega}(zg)
    \]
    pointwise. Moreover, according to  \cref{th-cv-psi},
    \[
        \psi_{g_n} \xrightarrow[n\to\infty]{} \psi_g
    \]
    pointwise. Thus, $\psi_g = 2\dz P^{\Omega}(zg)$.
\end{proof}

\section{The \texorpdfstring{$L^2$}{L²} problem}\label{sec-l2}

The aim of this section is to show how to use the $H^2$-results in order to prove the $L^2$-results.

\subsection{The space \texorpdfstring{$A^2 + \overline{A^2}$}{}}

We recall that for a bounded domain $\Omega \subset \C$, $A^2(\Omega)$ denotes the Bergman space on $\Omega$, i.e., the set of holomorphic function on $\Omega$ that are in $L^2(\Omega)$.

To prove \Cref{th-NCH-sub-NCL}, we need the following proposition:
\begin{prop}\label{th-berg+antiberg}
    Let $\Omega$ be an open bounded subset of $\C$. Assume that $\partial\Omega$ consists of a finite number of closed curves having a continuous tangent except at a finite number of corners of angle in $(0,2\pi)$. Then, $A^2(\Omega) + \overline{A^2(\Omega)}= L^2(\Omega)\cap h(\Omega)$ is closed.%  =  \{f\in L^2(\Omega,\ \C),\ f\text{ is harmonic}\}$.
\end{prop}
Actually, it is known~\cite[Corollary 4.3]{PS00II} that $A^2(\Omega) + \overline{A^2(\Omega)} = \{f\in L^2(\Omega,\ \C),\ f\text{ is harmonic}\}$.
As observed in \cite[Section 4]{PS00II} the proposition above is a consequence of the following inequality of Friedrichs~\cite{Friedrichs37}: if $\Omega$ satisfies the condition of the proposition, then there is $0<\theta<1$ such that for every $f\in A^2(\Omega)$ with $\int_{\Omega}f\diff A=0$, we have
\[
    \left|\int_{\Omega} f(z)^2 \diff A(z) \right| \leq \theta \int_\Omega |f(z)|^2 \diff A(z).
\]
See also \cite{Shapiro81}, \cite{PS00} for more on Friedrich's inequality and alternative proofs. 

For convenience of the reader, let us deduce \cref{th-berg+antiberg} from the Friedrichs inequality, using some ingredients from~\cite[Corollary 4.4]{PS00II}.

\begin{proof}[Proof of \cref{th-berg+antiberg} thanks to Friedrichs' inequality]
Let $H \coloneqq \{f\in A^2(\Omega),\ \int_\Omega f \diff A = 0\}$. $H$ is a closed subspace of $L^2(\Omega)$. Let $\Pi\colon L^2(\Omega) \to H$ be the orthogonal projection on $H$, and let $F$ be the antilinear operator $f\in H\mapsto \Pi \overline f \in H$. 

\step{$\|F\|\leq \theta$, where $\theta$ is the constant of Friedrichs' inequality}
Thanks to the definition of $F$, we get that for every $f,g\in H$,
\[
    (f,Fg)_{L^2} = \int_\Omega f(z)g(z) \diff A(z).
\]
In particular, $(f,Fg)_{L^2} = (g,Ff)_{L^2}$. Thus, thanks to the semilinearity of the scalar product,
\[
(f,Fg)_{L^2} = \frac14\big( (f+g,F(f+g))_{L^2} - (f-g,F(f-g))_{L^2}\big).
\]
Hence, according to Friedrichs' inequality and the parallelogram identity,
\[
|(f,Fg)_{L^2}| \leq \frac14\big( \theta \|f+g\|^2 + \theta\|f-g\|^2\big) = \frac\theta 4(2\|f\|^2 + 2\|g\|^2) = \frac\theta2 (\|f\|^2 +\|g\|^2).
\]
Finally, applying this to $f=\|g\|_{L^2}\|Fg\|_{L^2}^{-1}Fg$,
\[
\|Fg\| \leq \theta \|g\|.
\]

\step{Conclusion}
Let $f_n,g_n\in A^2(\Omega)$ be such that $h_n \coloneqq f_n +\overline{g_n}$ converges in $L^2(\Omega)$, and let us prove that the limit $h$ is in $A^2(\Omega) + \overline{A^2(\Omega)}$. Since $\lambda_n\coloneqq \int_{\Omega} h_n \diff A$ converges, we can assume without loss of generality that $\int_\Omega f_n \diff A = \int_\Omega g_n\diff A = 0$ (subtract their average otherwise).

According to Cauchy-Schwarz inequality and step 1, we have for every $f,g\in H$
\begin{align*}
    \|f+\overline{g}\|_{L^2}^2
    &= \|f\|_{L^2}^2 + \|g\|_{L^2}^2 +2 \Real \int_\Omega f(z) g(z) \diff A(z)\\
    &=\|f\|_{L^2}^2 + \|g\|_{L^2}^2 +2\Real(Ff,g)_{L^2}\\
    & \geq \|f\|_{L^2}^2 + \|g\|_{L^2}^2 -2\theta \|f\|_{L^2}\|g\|_{L^2}\\
    &\geq (1-\theta)(\|f\|_{L^2}^2 + \|g\|_{L^2}^2).
\end{align*}
Thus, $f_n$ and $g_n$ are Cauchy sequences, hence they converge in $A^2(\Omega)$. This proves that $h= \lim(f_n+\overline{g_n}) = \lim f_n+ \lim \overline{g_n} \in A^2(\Omega) + \overline{A^2(\Omega)}$.
\end{proof}

\Cref{th-berg+antiberg} allows us to prove the following result, where $L^2_0([0,T]\times \omega)$ is the space $\{f\in L^2([0,T]\times\omega),\ \int_{[0,T]\times \omega} f(t,\eu^{\iu x})\diff t\diff x = 0\}$. Recall the input-to-output map
for the $L^2$ control system~\eqref{eq-half-L2} defined on $L^2([0,T]\times \omega)$: %is $\FTp\coloneqq P_+\FT\colon L^2([0,T]\times\omega)\to H^2$, i.e.,
\[
    \FT u\coloneqq\int_0^{T}S(T-s)(\chi_{\omega}u(s,\cdot))\diff s,
\]
and correspondingly for \cref{eq-half-H2}:
\[
    \FTp=P_+\FT.
\]
\begin{corollary}\label{th-control-0-0}
    Let $T>0$ and
    %\footnote{The map $\FT$ is the input to output map of \Cref{def-Ft}.} 
    $\mathcal U \coloneqq \{ v\in L^2([0,T]\times \omega),\ \FTp v = 0 \}$ the set of controls $v$ that steer $0$ to $0$ in time $T$ for the system $(\partial_t + |D|)h = P_+ \mathds 1_\omega v$. We have
    \[\mathcal U + \overline{\mathcal U} = L^2_0([0,T]\times \omega).\]
\end{corollary}

\begin{proof}
    We have $\mathcal U = \ker(P_+ \FT) = \range((\FT)^*P_+)^\bot$. As already mentioned in \eqref{adjoint}, $(\FT)^*g(t) = \mathds 1_\omega^*S(T-t)^*g$,
    where $\mathds 1_\omega^*\colon L^2(\T) \to L^2(\omega)$ is the restriction operator. Hence, if $g_0\in L^2(\T)$ and $g_0^+ = P_+g_0 \in H^2$, $(\FT)^* P_+g_0(t,\eu^{\iu x}) = g_0^+(\eu^{-(T-t)+\iu x})$ for $0<t<T$ and $\eu^{\iu x} \in \omega$. Thus, setting $z=\eu^{-(T-t)+\iu x}$,
    \begin{align*}
        u\in \mathcal U
        & \eq \forall g_0^+ \in H^2, \int_{[0,T]\times\omega} u(t,\eu^{\iu x}) \overline{g_0^+(\eu^{-(T-t)+\iu x})} \diff t\diff x =0 \\
        & \eq \forall g_0^+ \in H^2, \int_{\Omega^*} |z|^{-2}u(-{\ln\lvert z\rvert}, z/|z|) \overline{g_0^+(z)} \diff A(z) = 0          \\
        & \eq u_\an \in (H^2)^\bot \text{ (orthogonal in $L^2({\Omega^*})$),}
        \intertext{where $u_\an(z) = |z|^{-2} u(-{\ln\lvert z\rvert}, z/|z|)$, and where as in the rest of the article, $\Omega^*= \Omega^*_T = \{z\in \C,\ \eu^{-T}<|z|<1,\ z/|z| \in \omega\}$. Since polynomials are dense in $H^2$ and in $A^2(\Omega^*)$~\cite[Chapter 18, Theorem 1.11]{Conway95}, we get}
        u\in \mathcal U & \eq u_\an \in (A^2(\Omega^*))^\bot.
    \end{align*}
    Hence, using the same notation $u_\an$ as above,
    \begin{align*}
        u \in \mathcal U + \overline{\mathcal U}
        & \eq u_\an \in A^2({\Omega^*})^\bot + \overline{A^2({\Omega}^*)}^\bot.
        \intertext{Since $A^2({\Omega^*}) + \overline{A^2({\Omega^*})}$ is closed (\Cref{th-berg+antiberg}), we get $A^2({\Omega^*})^\bot + \overline{A^2({\Omega^*})}^\bot = \big(A^2({\Omega^*}) \cap \overline{A^2({\Omega^*})}\big)^\bot$ \cite[Theorem~2.16]{Brezis11}. Hence, $A^2({\Omega^*})^\bot + \overline{A^2({\Omega^*})}^\bot = \C^\bot$ (orthogonal of constant functions in $L^2(\Omega^*)$). Hence,}
        u \in \mathcal U + \overline{\mathcal U} & \eq u_\an \in \C^\bot \\
        & \eq \int_{\Omega^*}|z|^{-2} u(-{\ln\lvert z\rvert},\arg(z/|z|)) \diff A(z) = 0 \\
        & \eq \int_{[0,T]\times \omega} u(t,\eu^{\iu x}) \diff t\diff x = 0.\qedhere
    \end{align*}
\end{proof}

\subsection{Link between the \texorpdfstring{$H^2$}{H²} and \texorpdfstring{$L^2$}{L²} control systems}
Before proving \cref{th-NCH-sub-NCL}, let us point out that
\begin{equation*}
    P_+\overline g(\eu^{\iu x}) 
    = P_+ \left(\sum_{n\in\Z} \overline{\hat{g}(n)} \eu^{-\iu nx}\right)
    = \sum_{n\leq 0} \overline{\widehat{g(n)}} \eu^{-\iu n x}.
\end{equation*}
Hence, as $P_+g$ is the projection on nonnegative frequencies of $g$, $\overline{P_+\overline g}$ is the projection on nonpositive frequencies, a fact we will use frequently in the rest of this section.

\begin{delayedproof}{th-NCH-sub-NCL}
    \step{The easy implication}
    If
    \[
        (\partial_t +|D|)f(t,\eu^{\iu x}) = \mathds 1_\omega u(t,\eu^{\iu x}),\qquad  f(0,\eu^{\iu x}) = f_0(\eu^{\iu x}),
    \]
    then
    \begin{equation*}%\label{eq-u_pm_control}
        \begin{array}{l}
            (\partial_t +|D|)P_+f(t,\eu^{\iu x}) = P_+ \mathds 1_\omega u(t,\eu^{\iu x}),\qquad  P_+f(0,\eu^{\iu x}) = P_+f_0(\eu^{\iu x}), \\
            (\partial_t +|D|)P_+\overline{f(t,\eu^{\iu x})} = P_+ \mathds 1_\omega \overline{u(t,\eu^{\iu x})},\qquad  P_+\overline{f(0,\eu^{\iu x})} = P_+\overline{f_0}(\eu^{\iu x}).
        \end{array}
    \end{equation*}

    Hence, if $f_0\in \NCL(\omega,T)$, then $P^+f_0 \in \NCH(\omega)$ and $P^+\overline{f_0}\in \NCH(\omega)$.

    \step{The less easy implication}
    Assume that $P_+f_0\in \NCH(\omega)$ and $P_+\overline{f_0}\in \NCH(\omega)$. Then there exists $u^\pm\in L^2([0,T]\times \omega)$ such that the solutions $f^\pm$ of
    \begin{equation}\label{eq-u_pm_control}
        \begin{array}{l}
            (\partial_t +|D|)f^+(t,\eu^{\iu x}) = P_+ \mathds 1_\omega u^+(t,\eu^{\iu x}),\qquad  f^+(0,\eu^{\iu x}) = P_+f_0(\eu^{\iu x}), \\
            (\partial_t +|D|)f^-(t,\eu^{\iu x}) = P_+ \mathds 1_\omega \overline{u^-(t,\eu^{\iu x})},\qquad  f^-(0,\eu^{\iu x}) = P_+\overline{f_0}(\eu^{\iu x}).
        \end{array}
    \end{equation}
    satisfy $f^\pm(T,\cdot) = 0$.

    We are looking for a $v$ that steers $f_0$ to $0$ for the $L^2$-control system~\eqref{eq-half-L2}. Let us denote by $f = S(t)f_0 + \FT[t] v$ the corresponding solution. Considering $P_+f$ and $\overline{P_+ \overline f}$, we see that $v$ steers $f_0$ to $0$ if and only if $v$ steers both
    \begin{itemize}
        \item $P_+f_0$ (the nonnegative frequencies of $f_0$) to $0$ for the control system $(\partial_t + |D|)h^+ = P_+ \mathds 1_\omega v$,
        \item $\overline{P_+\overline{f_0}}$ (i.e., the nonpositive frequencies of $f_0$) to $0$ for the control system $(\partial_t + |D|)h^- = \overline{P_+ \mathds 1_\omega \overline v}$.
    \end{itemize}
    
    The first of the above two conditions is equivalent to $v \in u^+ + \mathcal U$, while the second condition is equivalent to $v\in u^- +\overline{\mathcal U}$. Hence, we want to prove that $(u^++\mathcal U)\cap(u^- +\overline{\mathcal U}) \neq \emptyset$. This is equivalent to $u^+ - u^- \in \mathcal U + \overline{\mathcal U}$.

    But $\mathcal U + \overline{\mathcal U} = L^2_0([0,T]\times \omega)$ (\Cref{th-control-0-0}). Hence, we only have to prove that
    \begin{equation}\label{eq-NCL-NCH-subgoal}
        \int_{[0,T]\times \omega} u^+(t,\eu^{\iu x}) \diff t\diff x =
        \int_{[0,T]\times \omega} u^-(t,\eu^{\iu x}) \diff t\diff x
    \end{equation}

    Considering the $0$th Fourier coefficient in the control systems~\eqref{eq-u_pm_control}, we see that
    \[
        \partial_t \widehat{f^+}(t,0) = \widehat{\mathds 1_\omega u^+}(t,0).
    \]
    Hence,
    \begin{equation*}
        \int_{[0,T]\times\omega} u^+(t,\eu^{\iu x}) \diff t\diff x
        =\int_0^T \widehat{\mathds 1_\omega u^+}(t,0) \diff t=\widehat{f^+}(T,0) - \widehat{f^+}(0,0) = -\widehat{f_0}(0).
    \end{equation*}
    Similarly,
    \[
        \int_{[0,T]\times\omega} \overline{u^-}(t,\eu^{\iu x}) \diff t\diff x = - \widehat{\overline{f_0}}(0).
    \]
    Thus, we do have \cref{eq-NCL-NCH-subgoal}, and the proof of \Cref{th-NCH-sub-NCL} is complete.
\end{delayedproof}

\subsection{Properties of the \texorpdfstring{$L^2$}{} control system}
We now combine the properties of the $H^2$ control system and the link between the $H^2$ and $L^2$ control system (\cref{th-NCH-sub-NCL}).

\begin{delayedproof}{th-NCL-non-NC}
    Assume that $f_0\in \NCL(\omega)$ is analytic in a neighborhood of $\T$, in other words,
    \[
        |\widehat{f_0}(n)|\leq C\eu^{-c|n|}.
    \]

    In particular, $P_+f_0(z) = \sum_{n\geq 0} \widehat{f_0}(n) z^n$ extends holomorphically to $D(0,\eu^c)$. But $P_+f_0\in\NCH(\omega)$ (\cref{th-NCH-sub-NCL}). Hence, according to \cref{th-analy-non-nc} $P_+ f_0 = 0$.

    Similarly, $P_+ \overline{f_0} = \sum_{n\leq 0} \overline{\widehat{f_0}(n)} z^{-n}$ is in $\NCH(\omega)$ and extends analytically to $D(0,\eu^c)$, hence $P_+\overline{f_0} = 0$.

    Thus, Fourier coefficients for both nonnegative and nonpositive $n$ are $0$, and so $f_0=0$.
\end{delayedproof}

To prove that $\NCL$ is dense, we will often use that $P^+$ is self-adjoint on $W^{k,2}(\T)$. We will also use the following variation of \cref{th-lemma-ortho}:
\begin{lemma}\label{th-lemma-ortho-const}
let $k\in \N$. Let $f\in H^2\cap W^{k,2}(\T)$. Assume that
\begin{equation*}
    \forall g\in C_c^\infty(\omega)\text{ such that } \int_\omega g = 0,\ \langle f,g\rangle_{W^{k,2}(\T)}=0.
\end{equation*}
Then $f$ is constant.
\end{lemma}
\begin{proof}
Choose $\chi\in C_c^\infty(\omega)$ real valued with $\int_\omega \chi = 1$. Set $f_1 = f-\langle f,\chi \rangle_{W^{k,2}}$.

Let $g\in C_c^\infty(\omega)$. Set $g_1 = g - \chi\int_\omega g$, which satisfies $g_1\in C_c^\infty(\omega)$ and $\int_\omega g_1 = 0$.
\begin{align*}
    \langle f_1,g\rangle_{W^{k,2}}
    &= \langle f,g\rangle_{W^{k,2}} - \langle f,\chi \rangle_{W^{k,2}} \langle 1,g\rangle_{W^{k,2}}\\
    &= \langle f,g\rangle_{W^{k,2}} - \langle f,\langle g,1\rangle_{W^{k,2}}\chi \rangle_{W^{k,2}}\\
    &= \big\langle f,g-\langle g,1\rangle_{W^{k,2}}\chi\big\rangle_{W^{k,2}}.
    \intertext{Notice that according to the definition of the $W^{k,2}$ scalar product, $\int_\omega g = \langle g,1\rangle_{W^{k,2}}$. Thus, according to the hypothesis on $f$,}
    \langle f_1,g\rangle_{W^{k,2}} 
    &= \langle f,g_1\rangle_{W^{k,2}}\\
    &= 0.
\end{align*}
Since this is true for every $g\in C_c^\infty(\omega)$ and since $f_1\in H^2$, we can use \cref{th-lemma-ortho} to conclude that $f_1=0$, i.e., $f=\langle f,\chi\rangle_{W^{k,2}}$ is constant.
\end{proof}
\begin{delayedproof}{th-NCL-dense}
    Let $h\in \NCL(\omega)^\bot$ (orthogonal in $W^{k,2}(\T)$).
    
    \step{$P_+h$ is constant} Let $g\in C_c^\infty(\omega)$ with $\int_\T g = 0$ and set $g^+ = P_+ g$. According to \cref{th-CondSuff}, $P_+g^+ = P_+ g\in\NCH(\omega)$. Moreover, $P_+ \overline{g^+}(z) = \overline{\widehat g}(0)$ ($P_+\overline{g^+}$ only sees the nonpositive frequencies of $g^+$, but $g^+ = P_+ g$ does not have any negative frequencies). But we assumed that $\int_\T g = 0$, hence $P_+ \overline{g^+} = 0\in\NCH(\omega)$. Thus, $g^+ \in \NCL(\omega)$ (\cref{th-NCH-sub-NCL}).

    Thus, for every $g\in C_c^\infty(\omega)$ with $\int_\T g = 0$,
    \begin{equation*}
        0= \langle h,g^+\rangle_{W^{k,2}(\T)}
        =\langle h,P_+ g\rangle_{W^{k,2}(\T)}
        =\langle P_+ h, g\rangle_{W^{k,2}(\T)}
    \end{equation*}
    In addition, since $h\in L^2(\T)$, $P_+ h\in H^2$. Thus, according to \cref{th-lemma-ortho-const}, $P_+ h$ is constant.

    \step{$\overline{P_+\overline h}$ is constant} The reasoning is similar. Let $g\in C_c^\infty(\omega)$ with $\int_\T g = 0$, and set $g^- = \overline{P_+ \overline g}$. According to \cref{th-CondSuff}, $P_+ \overline{g^-} = P_+ \overline g\in\NCH(\omega)$. Moreover, $P_+ g^- (z) = \widehat g(0) = 0 \in\NCH(\omega)$. Thus, $g^- \in \NCL(\omega)$ (\cref{th-NCH-sub-NCL}).
    
    Thus, for every $g\in C_c^\infty(\omega)$ with $\int_\T g = 0$,
    \begin{equation*}
        0= \langle h, g^-\rangle_{W^{k,2}(\T)}
        =\langle h,\overline{P_+ \overline g}\rangle_{W^{k,2}(\T)}
        =\overline{\langle \overline h, P_+ \overline g\rangle_{W^{k,2}(\T)}} 
        = \overline{\langle P_+\overline h, \overline g\rangle_{W^{k,2}(\T)}}
        =\langle \overline{P_+\overline h}, g\rangle_{W^{k,2}(\T)}.
    \end{equation*}

    Thus $\overline{P_+ \overline h}$ is constant.

    \step{Conclusion} From the two first steps, every positive frequency of $h$ is $0$, and from the second step, every negative frequency is $0$, so that $h$ is constant.

    Finally, consider $g\in C_c^\infty(\omega)$ with $\int_\T g \neq 0$. According to \cref{th-CondSuff}, $P_+g\in \NCH(\omega)$ and $P_+ \overline g\in\NCH(\omega)$. Hence, according to \cref{th-NCH-sub-NCL}, $g\in \NCL(\omega)$. Consequently, $0 = \langle h,g\rangle_{W^{k,2}(\T)}$, and since $h$ is constant, $0 = \langle h,g\rangle_{W^{k,2}(\T)}=\langle h,g\rangle_{L^2}=\hat{h}(0)\int_{\T}g$, so that finally $h = 0$.

    We proved that $\NCL(\omega)^\bot = \{0\}$.
\end{delayedproof}

\appendix

\section{Construction of the cutoff function}\label{sec-cutoff}

\begin{lemma}
Let $\Omega_1$ and $\Omega_2$ be open subsets of $\C$ such that $\distance(\Omega_1\setminus \Omega_2\ \Omega_2\setminus \Omega_1) > 0$.
There exists a function $\chi \in C^\infty(\C)$ such that $\chi = 0$ on a neighborhood of $\cl{(\Omega_1\setminus \Omega_2)}$ and $\chi=1$ on a neighborhood of $\cl{(\Omega_2 \setminus \Omega_1)}$, with additionally $0\leq \chi \leq 1$ and
\begin{equation}
   \forall x\in \Omega_1\cup\Omega_2), |\nabla \chi(x)| \leq \dfrac{C}{\distance(\Omega_1\setminus \Omega_2,\ \Omega_2\setminus \Omega_1)}
\end{equation}
on $\Omega_1\cap \Omega_2$.
\end{lemma}

\begin{proof}
Set $\epsilon = \distance(\Omega_1\setminus \Omega_2,\ \Omega_2\setminus \Omega_1)$. Let $\rho \in C_c^\infty(B(0,1))$ with $\rho\geq 0$ and $\int_\C \rho(z) \diff A(z) = 1$. For $\eta >0$, set $\rho_\eta(z) = \eta^{-2} \rho(\eta^{-1}z)$.

Set $d_1(z) = \distance(z,\ \Omega_1\setminus \Omega_2)$ and $d_2(z) = \distance(z,\ \Omega_2\setminus \Omega_1)$. Set $\delta_i^\eta = d_i \ast \rho_\eta$.

Let $\varphi \in C^\infty([0,1])$ with $\varphi = 0$ on $[0,1/3]$ and $\varphi = 1$ on $[2/3,1]$.

Finally, for some $\eta >0$ to be chosen later (depending on $\epsilon$), set
\begin{equation}\label{eq-def-cutoff}
    \chi(z) = \varphi\left(\dfrac{\delta_1^\eta(z)}{\delta_1^\eta(z)+\delta_2^\eta(z)}\right).
\end{equation}

First, we claim that
\begin{equation}\label{eq-d-delta}
    \|d_i-\delta_i^\eta\|_{L^\infty} \leq C_1\eta.
\end{equation}
This a consequence of straightforward computations using that $d_i$ is $1$-lipschitz:
\begin{align*}
    | \delta_i^\eta(z) - d_i(z) |
    &\leq \int_\C| d_i(w) - d_i(z)|\rho_\eta(z-w) \diff A(w)\\
    &\leq \int_\C |z-w|\rho_\eta(z-w) \diff A(w)\\
    &= \eta \int_\C |w| \rho(w) \diff A(w).
\end{align*}
From now on, we choose $\eta = \epsilon/(16 C_1)$ (with the $C_1$ of the previous estimate~\eqref{eq-d-delta}), so that $\|d_i-\delta_i^\eta\|_{L^\infty} \leq  \epsilon/16$.

Here are some facts about these functions:
\begin{enumerate}
    \item $|\nabla d_i| \leq 1$, hence $
    \|\nabla \delta_i^\eta\|_{L^\infty} = \|\nabla d_i \ast \rho_\eta\|_{L^\infty}
    \leq \|\nabla d_i\|_{L^\infty} \|\rho_\eta\|_{L^1} \leq 1$.
    \item By definition of $\epsilon$, we cannot have both $d_1< \epsilon/2$ and $d_2< \epsilon/2$, hence $d_1+d_2 \geq \epsilon/2$, and finally $\delta_1^\eta + \delta_2^\eta \geq d_1+d_2 - \epsilon/8 \geq 3\epsilon/8$.
    \item If $d_1 \leq \epsilon/16$, we have $\delta_1^\eta \leq \epsilon/16+ \epsilon/16$, hence $\delta_1^\eta/(\delta_1^\eta + \delta_2^\eta) \leq \frac{\epsilon/8}{3\epsilon/8} = 1/3$. Hence, by definition of $\chi$ (\cref{eq-def-cutoff}), $\chi = 0$ on a neighborhood of $\cl{(\Omega_1\setminus \Omega_2)}$.
    \item Similarly, $\chi = 1$ on $\{z,\ d_2(z)<\epsilon/16\}$, which is a neighborhood of $\cl{(\Omega_2\setminus \Omega_1)}$.
\end{enumerate}

The only thing left to prove on $\chi$ is the gradient estimate~\eqref{eq-cutoff-est}. We have
\begin{align*}
    \partial_x \chi
    &= \partial_x \left(\frac{\delta_1^\eta}{\delta_1^\eta + \delta_2^\eta}\right) \varphi'\left(\frac{\delta_1^\eta}{\delta_1^\eta+\delta_2^\eta}\right)\\
    &=\frac{\delta_2^\eta\partial_x\delta_1^\eta - \delta_1^\eta\partial_x\delta_2^\eta}{(\delta_1^\eta + \delta_2^\eta)^2} \varphi'\left(\frac{\delta_1^\eta}{\delta_1^\eta+\delta_2^\eta}\right)\\
    &=\left(\frac{\partial_x\delta_1^\eta}{1+\delta_1^\eta/\delta_2^\eta} - \frac{ \partial_x\delta_2^\eta}{1+\delta_2^\eta/\delta_1^\eta}\right) \frac1{\delta_1^\eta+\delta_2^\eta} \varphi'\left(\frac{\delta_1^\eta}{\delta_1^\eta+\delta_2^\eta}\right).
\end{align*}
In $\Omega_1 \cap \Omega_2$, $\delta_i^\eta > 0$, hence $1+\delta_1^\eta/\delta_2^\eta \geq 1$. Hence, we can estimate the right-hand side as
\begin{align*}
    |\partial_x \chi|
    &\leq \left( |\partial_x \delta_1^\eta| + |\partial_x \delta_1^\eta|\right)
    \frac1{|\delta_1^\eta+\delta_2^\eta|} \|\varphi'\|_{L^\infty}.\\
    \intertext{Using that $|\nabla \delta_i^\eta|\leq 1$ and $\delta_1^\eta + \delta_2^\eta\geq 3\epsilon/8$, we get}
    |\partial_x \chi|
    &\leq \frac{16}{3\epsilon} \|\varphi'\|_{L^\infty}.
\end{align*}
The same holds for $\partial_y \chi$.
\end{proof}

\section{Background on the boundary behaviour of conformal maps and on the Poisson kernel}\label{sec-conformal-and-poisson}

In this section, $\Omega\subset \C$ satisfies:
\begin{equation}\tag{H}\label{hyp-Omega}
\begin{array}{l}
    \Omega \text{ is a Jordan domain, i.e., the bounded component of $\C\setminus \gamma$, where $\gamma$ is a closed Jordan curve,}\\
    \partial\Omega \text{ is a finite union of smooth (say $C^2$) curves that meets at angle $\pi \alpha_k$ with $0<\alpha_k<1$.}
\end{array}
\end{equation}
The cases we are interested in are $\Omega=\{\eu^{x}\zeta,\ 0<x<T,\ \zeta\in \omega\} = \Omega_T$. Let also $\varphi$ be a Riemann mapping from $\D$ to $\Omega$.

Based on some existing results about the boundary behaviour of $\varphi$, we state and prove the results we need about the Poisson kernel and the Bergman Kernel of $\Omega$.

We do not claim much novelty here: some of these results (or similar ones) are either in the literature in great generality or are probably folklore. But to help the reader, we still prove them in our specific case. Let us also mention that not every theorem that follows requires all the above hypotheses on $\Omega$, even with our proof strategy. Since our goal is not to get the most general statements possible, we do not detail this.

According to Carathéodory's theorem~\cite[Theorem 2.6]{Pommerenke92}, $\varphi$ extends by continuity to $\cl \D$, and this extension (still denoted $\varphi$) is a homeomorphism from $\cl \D$ to $\cl \Omega$. In fact, if $\zeta_0 = \varphi(z_0) \in \partial\Omega$ is a corner of angle $\pi \alpha$, then, we have the following asymptotics when $z\to z_0$~\cite[Theorem 3.9]{Pommerenke92}:
\begin{align}
    \varphi(z) - \varphi(z_0) & \sim c (z-z_0)^\alpha,\label{eq-est-corner}           \\
    \varphi'(z)               & \sim c'(z-z_0)^{\alpha-1}.\label{eq-est-corner-deriv}
\end{align}
Note that this applies also when $\partial\Omega$ is smooth at $\zeta_0$: this is the case when $\alpha = 1$.

Finally, let us mention that the linear measure $\lvert\diff z|$ on $\partial\Omega$ is given by $\lvert\diff z| = |\varphi'(w)|\lvert\diff w|$ \cite[Theorem 6.8]{Pommerenke92}, i.e., if $J\subset \partial \D$,
\begin{equation}\label{eq-linear-measure}
    \int_{\varphi(J)} \lvert\diff z| = \int_J |\varphi'(w)|\lvert\diff w|.
\end{equation}

Here are some useful properties on $\varphi$ that we will use in this appendix:
\begin{lemma}\label{th-varphi-lemma}
Assume $\Omega\subset\C$ satisfies the hypothesis~\eqref{hyp-Omega}.
    \begin{itemize}
        \item Let $\psi = \varphi^{-1}$. Then, $\psi'$ extends continuously to $\cl\Omega$.
        \item There exists $C>0$ such that for every $0<r<1$ and almost every $|w|=1$, $|\varphi'(rw)|\leq C |\varphi'(w)|$.
    \end{itemize}
\end{lemma}

\begin{proof}
    \step{First assertion}
    Of course, $\psi'$ is continuous and does not vanish inside $\Omega$. To prove the first point, we only need to look on a neighborhood of each point of the boundary. Let $\zeta_0 = \varphi(z_0) \in \partial\Omega$. Then, inverting the asymptotics~\eqref{eq-est-corner}, we get
    \[
        \psi(\zeta) - \psi(\zeta_0)\sim c^{-1/\alpha} (\zeta-\zeta_0)^{1/\alpha}.
    \]
    Plugging this asymptotic in the asymptotics of $\varphi'$, we get
    \begin{align*}
        \psi'(\zeta) & = \dfrac1{\varphi'(\psi(\zeta))}                                               \\
                     & \sim \dfrac1{c' \big(c^{-1/\alpha}(\zeta-\zeta_0)^{1/\alpha}\big)^{\alpha-1}} \\
                     & \sim c'' (\zeta-\zeta_0)^{1/\alpha-1}.
    \end{align*}
    Since $\alpha <1$ if $\zeta_0$ is a corner and $\alpha=1$ else, $\psi'$ has a limit in $\zeta_0$ (equal to zero at the corners). Since this is true for any $\zeta_0\in \partial\Omega$, $\psi'$ does extend continuously to $\cl\Omega$.
    % Since $\psi' = 1/\varphi'\circ\psi$, we only have to prove that $|\varphi'|$ is bounded from below. Of course, $\varphi'$ is continuous and does not vanish inside $\D$. To conclude, we only have to look near the boundary.

    % Let $\zeta_0 = \varphi(z_0) \in \partial\Omega$. Let us denote the angle of the boundary at $\zeta_0$ by $\pi\alpha$.  Our hypothesis on $\Omega$ is that $0<\alpha\leq 1$. In this case, the above asymptotics tells us that $|\varphi'|$ is bounded from below on a neighborhood of $\zeta_0$ in $\cl\Omega$. Since this is true for any $\zeta_0\in \partial\Omega$, $|\varphi'|$ is indeed bounded from below.

    \step{Second assertion}
    According to the previous step, $\varphi'$ is bounded from below. Hence for any $0<r<1$, we have on $D(0,r)$, for every $|z|=1$, $|\varphi'(rz)| \leq C_r |\varphi'(z)|$.

    Let $|z_0|=1$. If $z\to z_0$, we have again $\varphi'(z) \sim c(z-z_0)^{\alpha-1}$ for some $0<\alpha\leq1$. In particular, there exists a neighborhood $U_{z_0}$ of $z_0$ in $\overline \D$, and $C_{z_0}>0$ such that
    \[
        \forall z \in U_{z_0},\ C_{z_0}^{-1} \varphi'(z) \leq |z-z_0|^{\alpha} \leq C_{z_0} \varphi'(z).
    \]
    
    With $|w| = 1$ and $0\leq r \leq 1$, we have
    \begin{align*}
        |w-z_0| \leq |w-rw| + |rw-z_0| \leq 2 |rw-z_0|.
    \end{align*}
    Since $0<\alpha \leq 1$, this proves
    \[
        |rw-z_0|^{\alpha -1} \leq 2^{1-\alpha}|w-z_0|^{\alpha-1}.
    \]
    
    Since $\varphi'(rw)$ is upper and lower-bounded by such a function on $U_{z_0}\cap \D$, there exists $C>0$ such that for every $rw\in U_{z_0}$
    \[
        |\varphi'(rw)| 
        \leq C_{z_0}|rw-z_0|^{\alpha-1} 
        \leq 2^{1-\alpha}C_{z_0}|w-z_0|^{\alpha-1}
        \leq 2^{1-\alpha}C_{z_0}^2 |\varphi'(w)|.
    \]
    Choosing a cover of $\T$ by a finite number of those $U_{z_0}$ yields the desired result.
\end{proof}

As a first consequence, here is a property on the Bergman kernel that is used in the proof of our necessary and sufficient condition (\cref{th-CNS})
\begin{prop}\label{th-berg-continuous}
Assume $\Omega\subset\C$ satisfies the hypothesis~\eqref{hyp-Omega}.
    Let $w\in \Omega$ and let $k^\Omega_w$ be the Bergman kernel (i.e., for all $f\in A^2(\Omega)$, $f(w) = \langle f,k_w\rangle_{A^2(\Omega)}$). Then, for all $w\in\Omega$, $k_w^\Omega$ extends continuously on $\overline\Omega$.
\end{prop}

\begin{proof}
    Let $k^\D_w(z) = (\pi(1-z \overline w)^2)^{-1}$ be the Bergman kernel of the unit disk~\cites[Theorem 2.7]{Zhu05}[\S1.2]{DS04}. Set $\psi = \varphi^{-1}$. By the change of variables formula \eqref{ChangeVarKernel}, 
    \[
        k^\Omega_w(z) = k_{\psi(w)}^\D(\psi(z)){\psi'(z)}\overline{\psi'(w)}.
    \]
    According to the formula for $k^\D$, for every $w \in \D$, $k_w^\D$ is continuous on $\cl\D$. Since $\psi$ is continuous on $\cl \D$ (Carathéordory's theorem) and $\psi'$ is also continuous on $\cl \D$ (\cref{th-varphi-lemma}), $k^\Omega_w$ is indeed continuous on $\cl\D$.
\end{proof}

% If $z\in \C$, we denote the $\diff m_z$ the harmonic measure on $\partial\Omega$. That is, if $g\in C^0(\partial\Omega)$, the generalised solution $u$ of
% \[
% \left\{ \begin{array}{ll}
%      -\Delta u(z) = 0,& z\in \Omega\\
%      u(w) = g(w),& w\in \partial_\Omega
% \end{array}\right.
% \]
% is
% \[
% u(z) = \int_{\partial_\Omega} g(w) \diff m_z(w).
% \]
% We denote this $u$ by $P^\Omega g$.

% \begin{theorem}
%  Let $\varphi\colon \D \to \Omega$ be a Riemann mapping. Then, $\diff m_z$ is given by 
%  \[\diff m_z(x) = P^\D(\varphi^{-1}(x),\varphi^{-1}(z))|(\varphi^{-1})'(x)|\diff \sigma(x).\]
%  Here, $P^\D$ is the Poisson kernel of $\D$.
% \end{theorem}

Another consequence is that the Dirichlet problem is solvable on $\Omega$ with $L^2$ boundary data. Recall that a function $h$ is harmonic in $\Omega$ if and only if $h\circ \varphi$ is harmonic in $\D$. In particular, if $h$ is regular enough (say $C^2(\Omega)\cap C^0(\overline \Omega)$), $g \coloneqq h_{|\partial\Omega}$ and $\Delta h = 0$, then, $\widetilde h = h\circ\varphi$ satisfies $\widetilde h\in C^2(\D)\cap C^0(\overline \D)$ (assuming hypothesis \eqref{hyp-Omega}), $\Delta \widetilde h = 0$ and $\widetilde h_{|\partial\D} = h\circ \varphi$. I.e., $h\circ \varphi = P^\D(g\circ\varphi)$. This motivates the following definition of $P^\Omega$:
\begin{corollary}\label{th-poisson-lipschitz}
Assume $\Omega\subset\C$ satisfies the hypothesis~\eqref{hyp-Omega}.
    Let $P^\Omega$ be the map $C^0(\partial\Omega)\to C^0(\cl\Omega)$ defined by
    \[
        P^\Omega g (\varphi(w)) = P^\D(g\circ\varphi)(w),\quad w\in\D,
    \]
    where $P^\D$ is the Poisson kernel of the disk:
    \[
     P^{\D}(z,\zeta)=\frac{1-|z|^2}{|z-\zeta|^2},\quad z\in\D,\zeta\in\T.
    \]

    The formula above extends $P^\Omega$ as a continuous linear map $L^2(\partial\Omega) \to h(\Omega)$, where $h(\Omega)$ is the space of harmonic functions on $\Omega$ endowed with the topology of uniform convergence on compact subsets.
\end{corollary}

We have already encountered the Poisson kernel for domains other than $\D$ when discussing harmonic extension to $\Dext$, but this domain does not satisfy the hypotheses~\eqref{hyp-Omega}.

By definition, we say that $u=P^\Omega g$ is the solution of $\Delta u = 0$ on $\Omega$, $u=g$ on $\partial\Omega$. The Dirichlet problem in nonsmooth domains has been widely studied, with hypotheses much more general than~\eqref{hyp-Omega}, see, e.g.,~\cite{Kenig94} for an overview. But in our case, we can prove what we need very easily.
\begin{proof}
    As in \cref{th-varphi-lemma}, we denote $\psi = \varphi^{-1}$. The definition of $P^\Omega$ reads
    \begin{IEEEeqnarray*}{rCl}
        P^\Omega g(\varphi(w))
        &=& \int_{\partial\D} P^\D(w,w') g(\varphi(w')) \lvert\diff w'|, \quad w\in \D.
        \intertext{According to the change of variable-formula for $\lvert \diff w'|$~\eqref{eq-linear-measure}, this can be written as}
        P^\Omega g(z)
        &=& \int_{\partial\Omega} P^\D (\psi(z),\psi(z')) g(z') |\psi'(z')|\lvert\diff z'|, \quad z=\varphi(w) \in \Omega.
    \end{IEEEeqnarray*}
    Since $|\psi'|$ is bounded~(\cref{th-varphi-lemma}), the definition of $P^\Omega g$ makes sense for $g\in L^2(\partial\Omega)$, and for every such $g$, $P^\Omega g$ is harmonic on $\Omega$. The required continuity property of the map $g\mapsto P^\Omega g$ follows from the fact that $P^\D(\psi(z),\psi(z')|\psi'(z')|$ is bounded on $D(0,r)_z\times \partial \D_{z'}$ for every $0<r<1$.
\end{proof}

The following proposition is a simple elliptic regularity result. Since it is enough for our purpose and since the proof is simple with the properties of the  conformal map outlined above, we do it here.
\begin{prop}\label{th-elliptic-reg}
Assume $\Omega\subset\C$ satisfies the hypothesis~\eqref{hyp-Omega}.
    Let $g\colon \partial\Omega\to \C$. For every $\zeta_0\in\partial\Omega$ such that $\partial\Omega$ is not smooth at $\zeta_0$, assume that $g=0$ on a neighborhood of $\zeta_0$. If $g\in W^{1,2}(\partial\Omega)$, then $P^\Omega g\in W^{1,2}(\Omega)$.
\end{prop}
\begin{proof}
    Let $\zeta_0 = \varphi(z_0) \in \partial\Omega$ that is not a corner of $\partial\Omega$. According to the asymptotics of $\varphi$ and $\varphi'$~\eqref{eq-est-corner}--\eqref{eq-est-corner-deriv}, $\varphi$ is $C^1$ at $z_0$.
    
    Since $g=0$ on a neighborhood of the corners, we deduce that $g\circ\varphi \in W^{1,2}(\T)$. We claim that $\dz P^\D(g\circ\varphi)\in L^2(\D)$. This is a standard computation: set $g_1 = g\circ\varphi$. We can write
    \[
        P^\D g_1(z)
        = \sum_{n\geq 0} \widehat{g_1}(n) z^n + \sum_{n>0} \widehat{g_1}(-n) \overline z^{n.}
    \]
    In particular,
    \[
     \|\dz (P^\D g_1)\|_{L^2(\D)}^2 
     = \bigg\|\sum_{n\geq 0} n\widehat{g_1}(n) z^{n-1}\bigg\|_{L^2(\D)}^2 
     = \pi\sum_{n> 0} n|\widehat{g_1}(n)|^2.
    \]
    Since $g_1\in W^{1,2}(\T)$, the right-hand side above is indeed finite (in fact, we only need $g_1\in W^{1/2,2}(\T)$). This proves the claim.
    
    Finally, the change of variables $z = \psi(\zeta) = \varphi^{-1}(\zeta)$ that satisfies $\diff A(z) = |\psi'(\zeta)|^2\diff A(\zeta)$ gives
    \begin{align*}
        \int_{\Omega}|\dz[\zeta] P^\Omega g(\zeta)|^2 \diff A(\zeta)
        &= \int_\Omega \left|\frac{\diff}{\diff \zeta}\left(P^\D g_1(\psi(\zeta))\right)\right|^2 \diff A(\zeta) \\
        &= \int_\Omega \left|\psi'(\zeta)\left(\dz P^\D g_1(\psi(\zeta))\right)\right|^2 \diff A(\zeta)\\
        &= \int_\D \left|\dz P^\D g_1(z)\right|^2 \diff A(z).
    \end{align*}
    This right-hand side is finite, hence $\dz[\zeta] P^\Omega g\in L^2(\Omega)$. Similarly, $\dzb[\zeta] P^\Omega g\in L^2(\Omega)$.
\end{proof}

\begin{prop}\label{th-weak-stokes}
Assume $\Omega\subset\C$ satisfies the hypothesis~\eqref{hyp-Omega}.
    We use the notation $P^\Omega$ of \cref{th-poisson-lipschitz}. Let $g\in L^2(\partial\Omega)$. Assume that $\dz P^\Omega g \in L^2(\Omega)$. Then, for every $v$ holomorphic on $\Omega$ that is $C^0(\cl\Omega)$, we have
    \[
        \int_{\partial\Omega} v(z)\overline g(z) \diff z = 2\iu\int_\Omega v(z)\overline {\dz P^\Omega g}(z) \diff A(z).
    \]
\end{prop}

If $\Omega$ is smooth and $f\in C^1(\cl\Omega)$, Stokes' theorem applied to the differential form $f(z)\diff z$, gives
\[
\int_{\partial\Omega} f(z)\diff z 
= \int_\Omega \dzb f(z) \diff\overline z\wedge \diff z
=2\iu \int_\Omega \dzb f(z) \diff A(z).
\]
If $v\overline{P^\Omega g} \in C^1(\cl\Omega)$, \cref{th-weak-stokes} is just this formula with $f = v \overline{P^\Omega g}$. The only difference in the statement above is the weaker regularity on $\Omega$ and $g$.

\begin{proof}
    We denote $u=P^\Omega g$. For $0<r<1$, let $\Omega_r = \varphi(D(0,r))$ that is a smooth domain. We apply the usual Stokes formula and take the limit $r\to 1$.

    \step{Stokes formula on $\Omega_r$}
    Since $\cl{(\Omega_r)}\subset \Omega$, $u \in C^1(\cl{(\Omega_r)})$. Moreover, $\Omega_r$ is smooth. Hence, according to Stokes' formula,
    \begin{equation}
        \int_{\partial\Omega_r} v(z) \overline u(z) \diff z 
        = 2\iu\int_{\Omega_r} \dzb(v\overline u)(z) \diff A(z)
        = 2\iu\int_{\Omega_r} v(z)\overline{\dz u}(z) \diff A(z), \label{eq-stokes-r}
    \end{equation}
    where the last equality follows from the holomorphy of $v$.

    \step{Limit of the right-hand side}
    Since $\dz u\in L^2(\Omega)$ (by hypothesis) and $v\in C^0(\cl{\Omega})$, this proves that $\mathds1_{\Omega_r}v\overline{\dz u}$ is dominated by an $L^1(\Omega)$ function. Hence, the dominated convergence theorem gives
    \[
        \int_{\Omega_r} v(z)\overline{\dz u}(z) \diff A(z) 
        = \int_{\Omega} \mathds1_{\Omega_r}v(z)\overline{\dz u}(z) \diff A(z)
        \xrightarrow[r\to 1]{} \int_{\Omega} v(z)\overline{\dz u}(z) \diff A(z).
    \]

    \step{Limit of the left-hand side}
    For $|z|<1$, let us set $U(z) = u(\varphi(z))$ and $V(z) = zv(\varphi(z)) \varphi'(z)$.

    We claim that
    \begin{equation}\label{eq-poisson-scal-prod}
        \int_{\partial{\Omega_r}} \overline u(z) v(z) \diff z = \iu \int_{0}^{2\pi} \overline g(\varphi(\eu^{\iu\theta})) v(\varphi(r^2\eu^{\iu\theta})) r^2 \eu^{\iu\theta} \varphi'(r^2\eu^{\iu\theta}) \diff \theta.
    \end{equation}

    Indeed, as already seen in the proof of \cref{th-varphi-lemma} the function $|\varphi'|$ is bounded from below by some constant $c>0$, so that
    \[
     \infty>\int_{\partial\Omega}|g(z)|^2\lvert\diff z|=\int_{\T}|g\circ\varphi(\zeta)|^2 |\varphi'(\zeta)| \lvert\diff \zeta|
     \ge c\int_{\T}|g\circ\varphi(\zeta)|^2  \lvert\diff \zeta|.
    \]
    Hence $g\circ\varphi\in L^2(\T)$. Then, since $U$ is harmonic in $\D$ with boundary values $g\circ\varphi$, and using $P^{\D}(z,\eu^{\iu\theta})$, the Poisson kernel on $\D$ already introduced earlier, we have
    \[
        U(z) = \int_0^{2\pi} g(\varphi(\eu^{\iu\theta})) P(z,\eu^{\iu\theta})\diff \theta.
    \]
    On the other hand
    \begin{align*}
        \int_{\partial\Omega_r} \overline u(z) v(z) \diff z
         & = \int_0^{2\pi} \overline u(\varphi(r\eu^{\iu t})) v(\varphi(r\eu^{\iu t})) r\iu \eu^{\iu t} \varphi'(r\eu^{\iu t}) \diff t                 \\
         & = \iu \int_{0}^{2\pi} \overline U(r\eu^{\iu t}) V(r\eu^{\iu t}) \diff t                                                                     \\
         & = \iu \int_0^{2\pi}\int_0^{2\pi} \overline g(\varphi(\eu^{\iu\theta})) P(r\eu^{\iu t},\eu^{\iu \theta})\diff \theta\ V(r\eu^{\iu t})\diff t .
        \intertext{Since $v$ is holomorphic, $V$ is harmonic. Using Fubini's theorem, we recognise a Poisson integral in the $t$~variable:}
        \int_{\partial\Omega_r} \overline u(z) v(z) \diff z
         & = \iu \int_0^{2\pi}\overline{g}(\varphi(\eu^{\iu\theta})) \int_0^{2\pi}   \underbrace{P(r\eu^{\iu t},\eu^{\iu \theta})}_{=P(r\eu^{\iu \theta},\eu^{\iu t})}V(r\eu^{\iu t})\diff t \diff \theta  \\
         & = \iu \int_0^{2\pi} \overline{g}(\varphi(\eu^{\iu\theta})) V(r^2 \eu^{\iu\theta}) \diff \theta,
    \end{align*}
    which is the claimed formula.

    We now justify that the right-hand side of this formula converges thanks to the dominated convergence theorem. According to \cref{th-varphi-lemma}
    \[
        |\varphi'(r\eu^{\iu \theta})|\leq C |\varphi'(\eu^{\iu \theta})|.
    \]

    The hypothesis $g\in L^2(\partial\Omega)$ reads $\int_0^{2\pi} |g(\varphi(\eu^{\iu \theta}))|^2 |\varphi'(\eu^{\iu \theta})|\diff \theta < +\infty$. Hence, the integrand of the right-hand side of \cref{eq-poisson-scal-prod} is dominated by the $L^1$ function $C |g(\varphi(\eu^{\iu\theta}))||\varphi'(\eu^{\iu \theta})|$, and the dominated convergence theorem gives
    \begin{align*}
        \lim_{r\to 1}\int_{\partial\Omega_r} \overline u(z) v(z) \diff z
         & = \iu \int_0^{2\pi} \overline g(\varphi(\eu^{\iu\theta})) v(\varphi(\eu^{\iu\theta})) \eu^{\iu\theta} \varphi'(\eu^{\iu\theta}) \diff \theta \\
         & = \int_{\partial\Omega} \overline g(z) v(z) \diff z.
    \end{align*}

    Combining this limit with the limit on $\int_{\Omega_r} \dzb (\overline uv)$ of the previous step, the theorem is proved.
\end{proof}

\section{The reachable space}
Some of our methods can be used to study the space of \emph{reachable states}, defined as follows:
\begin{definition}
    We say that $f_T\in H^2(\T)$ (resp.~$L^2(\T)$) is \emph{reachable on $\omega$ in time $T$ on $H^2$} (resp.~$L^2$) if there exists a control $u\in L^2([0,T]\times \omega)$ such that the solution $f$ of the control system~\eqref{eq-half-H2} (resp.~\eqref{eq-half-L2}) with initial condition $0$ satisfies $f(T,\cdot) = f_T$.

    We denote the set of reachable functions on $\omega$ in time $T$ by $\ReachH(\omega,T)$ (resp.~$\ReachL(\omega,T)$).
\end{definition}
With the methods of \cref{sec-l2}, we could prove results on $\ReachL(\omega,T)$ from results on $\ReachH(\omega,T)$. We do not detail this further, and instead sketch some results on $\ReachH(\omega,T)$.

\begin{prop}\label{th-MainEstim-reach}
    Let $\omega$ be an open subset of $\set T$ and $T>0$. Let $\Omega_T^* \coloneqq \{ z\in \set C,\ \eu^{-T}<|z|<1,\ z/|z| \in \omega\}$ (see \Cref{fig-Omega}). Let $f_T \in H^2$. The following assertions are equivalent:
    \begin{itemize}
        \item $f_T\in \ReachH(\omega,T)$.
        \item For every polynomial $g\in \C[X]$,
              \begin{equation}\label{MainEstim-Reach}
                  |\langle f_T,g\rangle_{H^2}| \lesssim \|g\|_{A^2(\Omega_T^*)}.
              \end{equation}
    \end{itemize}
\end{prop}
The proof is very similar to the one of \cref{MainEstim}, with the appropriate observability inequality.

\begin{theorem}\label{th-reach-holom-necessary}
    Let $T>0$ and $\omega$ be a strict open subset of $\T$. Let $\Omega_T$ as in \Cref{th-MainEstim}. If $f_T \in \ReachH(\omega)$, $f_T$ can be extended holomorphically on $\Cinf\setminus \cl{(\Omega_T)}$.
\end{theorem}
The proof is very similar to the one of \cref{th-NCH-holom-necessary} and is omitted.

We proved that every function in $P_+(W^{1/2,2}_0(\omega))$ is null-controllable to $0$. In fact, the same strategy based on Stokes formula used in \cref{sec-formal-stokes} proves these functions are also reachable:
\begin{theorem}\label{th-ReachH-lower}
    Let $T>0$ and $\omega$ be an interval of $\T$. If $h\in W^{1/2,2}_0(\omega)$, we implicitly extend $h$ on $\T$ by $0$ outside $\omega$. Then,
    \[
        P_+(W^{1/2,2}_{00}(\omega)) \subset \ReachH(\omega,T).
    \]
\end{theorem}

We might be able to adapt the methods used in the proof of \cref{th-CondSuff} to get a space of reachable states with less regularity, but we leave that for another time. However, we can exhibit explicit functions in $\ReachH(\omega,T)$ that are not null-controllable:
\begin{theorem}
    Let $T>0$ and $\omega$ be an open subset of $\T$. Let $\Omega_T^*$ as in \Cref{th-MainEstim-reach}. For $|u|<1$, let $k_u(z) = 1/(1-\overline{u}z)$. Then,
    \[
        k_u\in \ReachH(\omega,T) \eq u\in \Omega_T^*.
    \]
    In addition, $k_u\notin\NCH(\omega)$.
\end{theorem}
\begin{proof}
Since $k_u$ is the reproducing kernel of $H^2$, the left-hand side of the observability inequality~\eqref{MainEstim-Reach} is
\[
\langle k_u, g\rangle_{H^2} = g(u).
\]
This can be estimated by $\|g\|_{A^2(\Omega_T^*)}$ if and only if $u\in \Omega_T^*$.

Since $k_u$ is holomorphic on $D(0,1/|u|)$, it cannot be steered to $0$ (\cref{th-analy-non-nc}).
\end{proof}

In particular $\ReachH(\omega,T)$ genuinely depends on $T$. This is coherent with the following observation:
\begin{prop}
    \[
        \NCH(\omega) = S(T)^{-1}(\ReachH(\omega,T)).
    \]
\end{prop}
\begin{proof}
Let us denote by $S(t,f_0,u)$ the solution of the half-heat equation~\eqref{eq-half-H2}. According to Duhamel's formula, $S(t,f_0,u) = S(t)f_0 + \FTp[t] u$. Then, by definition of $\NCH(\omega)$ and $\ReachH(\omega,T)$,
\begin{align*}
    f_0\in \NCH(\omega)
    &\eq \exists u\in L^2([0,T]\times\omega),\ S(T,f_0,u) = 0\\
    &\eq \exists u\in L^2([0,T]\times\omega),\ S(T) f_0 = -\FTp u\\
    &\eq S(T) f_0 \in \ReachH(\omega,T).\qedhere
\end{align*}
\end{proof}

\subsection*{Acknowledgements} The authors express their gratitude to their colleagues for many interesting discussions on the topics, in particular Jérémi Dardé, Sylvain Ervedoza, Pierre Lissy and Pascal Thomas.
\printbibliography
\end{document}